\documentclass [11pt, twoside, leqno] {amsart}

\usepackage{amssymb, amsmath, amsfonts, amsthm, amsbsy, amscd}

\newtheorem{theorem}{Theorem}[section]
\newtheorem{proposition}{Proposition}[section]
\newtheorem{corollary}{Corollary}[section]
\newtheorem{lemma}{Lemma}[section]

\newtheorem{condition}{}{}

\theoremstyle{definition}
\newtheorem{example}{Example}[section]
\newtheorem{remark}{Remark}[section]

\numberwithin{equation}{section}

\newcommand{\Hom}{\text{Hom}}
\newcommand{\sy}{\mathcal{S}}
\newcommand{\cmzinf}{m^{(\infty)}}
\newcommand{\cmzmu}{m^{(\mu)}}
\newcommand{\cmzmub}{m^{(-2-\mu)}}
\newcommand{\re}{\mathfrak{Re}}

\newcommand{\pdop}{\mathcal{Q}}
\newcommand{\bih}{\Bar{\epsilon}}
\newcommand{\nat}{\mathbb{N}}
\newcommand{\cinfstar}{\,_{\infty}{\star}^{c}\,\,}
\newcommand{\tstar}{\,_{t}{\star}^{c/t}\,\,}

\newcommand{\hcmustar}{\,_{\epsilon}{\star}^{-\mu-n-1}\,\,}
\newcommand{\hmustar}{\,_{\epsilon}{\star}^{\mu}\,\,}
\newcommand{\chmustar}{\,_{\Bar{\epsilon}}{\star}^{-\mu-n-1}\,\,}
\newcommand{\chstar}{\,_{\Bar{\epsilon}}{\star}\,\,}
\newcommand{\bhmustar}{\,_{\Bar{\epsilon}}{\star}^{\mu}\,\,}
\newcommand{\hmubox}{\,_{\epsilon}{\square}^{\mu}\,\,}
\newcommand{\hmubbox}{\,_{\epsilon}{\square}^{-2-\mu}\,\,}
\newcommand{\hmubstar}{\,_{\epsilon}{\star}^{-2-\mu}\,\,}
\newcommand{\Sig}{\Sigma}

\newcommand{\si}{\sigma}
\newcommand{\tk}{\tilde{k}}
\newcommand{\T}{\mathcal{T}}
\newcommand{\ppois}{\langle \,,\, \rangle}
\newcommand{\la}{\langle}
\newcommand{\ra}{\rangle}

\newcommand{\hD}{\Hat{D}}
\newcommand{\hPhi}{\Hat{\Phi}}
\newcommand{\hteul}{\Hat{\mathbb{E}}}
\newcommand{\teul}{\mathbb{E}}
\newcommand{\heul}{\mathcal{X}}
\newcommand{\bnabla}{\bar{\nabla}}

\newcommand{\algh}{\mathcal{A}[[\epsilon]]}
\newcommand{\aff}{\mathcal{A}}
\newcommand{\ih}{\epsilon}
\newcommand{\thop}{\mathcal{L}_{t}}
\newcommand{\hop}{\mathcal{L}_{\epsilon}}
\newcommand{\hstar}{\,_{\epsilon}{\star}\,\,}
\newcommand{\hminusstar}{\,_{-\epsilon}{\star}\,\,}
\newcommand{\bhstar}{\,_{\epsilon}\Bar{\star}\,\,}
\newcommand{\pol}{\text{Pol}^{b}}
\newcommand{\polb}{\text{Pol}}
\newcommand{\dop}{\mathcal{D}}
\newcommand{\rc}{\mathcal{R}}
\newcommand{\op}{\mathcal{L}}
\newcommand{\kop}{\mathcal{K}}

\newcommand{\ga}{\gamma}
\newcommand{\al}{\alpha}
\newcommand{\be}{\beta}
\newcommand{\sklb}{\Bar{\mathcal{E}}^{(i_{1}\dots i_{k})}[\lambda]}
\newcommand{\skl}{\mathcal{E}^{(i_{1}\dots i_{k})}[\lambda]}
\newcommand{\en}{[\nabla]}
\newcommand{\ben}{[\Bar{\nabla}]}

\newcommand{\svol}{\Psi}

\newcommand{\standrepz}{\mathbb{V}_{0}}

\newcommand{\graded}{Gr}

\newcommand{\emf}{\mathcal{E}}
\newcommand{\emfb}{\Bar{\mathcal{E}}}

\newcommand{\hnabla}{\hat{\nabla}}

\newcommand{\standrep}{\mathbb{V}}

\newcommand{\form}{\mathsf{L}}
\newcommand{\dens}{\mathsf{L}}

\newcommand{\eul}{\mathbb{X}}

\newcommand{\tlop}{\mathcal{T}}
\newcommand{\proj}{\mathbb{P}}

\newcommand{\sdiv}{\text{div}}

\newcommand{\lie}{\mathfrak{L}}

\newcommand{\vect}{\text{Vec}}

\newcommand{\diff}{\text{Diff}}

\newcommand{\p}{\mathfrak{p}}

\newcommand{\tensor}{\otimes}

\newcommand{\rea}{\mathbb R}

\newcommand{\tr}{\text{tr} \,}

\newcommand{\aut}{\text{Aut}}

\newcommand{\pr}{\partial}
\newcommand{\diver}{\text{div}}
\newcommand{\alg}{\mathcal{A}} 
\newcommand{\balg}{\mathcal{B}}
\newcommand{\pois}{\{\,,\,\}}

\newcommand{\affdiff}{\mathcal{E}_{(ij)}^{k}}

\begin{document}
\title{Projectively Invariant Star Products}
\author{Daniel J. F. Fox} 
\email{fox@math.gatech.edu}
\address{School of Mathematics\\ Georgia Institute of Technology\\ 686 Cherry St.\\ Atlanta, GA 30332-0160, U.S.A}
% 53a55 is differential invariants 53a20 is projective differential geometry % 53C15 (general geometric structures)
\begin{abstract}
  It is shown that a (curved) projective structure on a smooth manifold determines on the Poisson algebra of smooth functions on the cotangent bundle, fiberwise-polynomial of bounded degree, a one-parameter family of graded star products. For a particular value of the parameter (corresponding to half-densities) the star product is symmetric, and specializes in the projectively flat case to the one constructed previously by C. Duval, P. Lecomte and V. Ovsienko. These star products are built from a projectively invariant quantization map associating to a symmetric polyvector a formally self-adjoint operator on densities. A limiting form of this family of star products yields a commutative deformation of the symmetric tensor algebra of the manifold which is closely related to the limiting commutative multiplication of modular forms defined by P. Cohen, Y. Manin, and D. Zagier. A basic ingredient of the proof is the construction of projectively invariant multilinear differential operators on bundles of weighted symmetric $k$-vectors. The construction works except for a discrete set of excluded weights and generalizes the Rankin-Cohen brackets of modular forms.
\end{abstract}
\thanks{During the preparation of this article the author benefited from a visit to the Instituto de Matem\'aticas y F\'isica Fundamental (IMAFF) del Consejo Superior de Investigaciones Cient\'ificas (CSIC) in Madrid, Spain.}

\maketitle

\section{Introduction}
The goal of this paper is to associate to (curved) projective structures on a smooth manifold, $M$, certain algebraic structures. Precisely, to each projective structure is associated a representation of symmetric polyvectors as self-adjoint operators on the space of half-densities on $M$, and from this is built a canonical deformation of the Poisson algebra structure on the algebra, $\pol(T^{\ast}M)$, of smooth functions on the cotangent bundle, polynomial in the fibers of globally bounded degree. These deformations have a close relation to the family of noncommutative multiplicative structures on the space of modular forms found by P. Cohen, Y. Manin and D. Zagier in \cite{Cohen-Manin-Zagier}, and generalize to the non-flat setting earlier work of C. Duval, P. Lecomte, V. Ovsienko and their collaborators. (That such a generalization should be possible was proposed in \cite{Duval-ElGradechi-Ovsienko}). The hope is that when properly understood these algebraic structures will point towards an algebraic generalization of the notion of projective structure. Because it works fairly well, the language used is that of deformation quantization, e.g. star products, quantization maps, and so on. It is not clear how natural is this language from a geometric point of view; in particular the notion of isomorphism usually used in deformation quantization is evidently too flabby (topological) to distinguish the star products associated to inequivalent projective structures.

In this introductory section the main theorems are stated. Precise definitions of some of the objects involved appear at the appropriate places in the main body of the text.

If $\alg$ is an associative algebra, let $\alg[[\ih]]$ be the ring of formal power series in the indeterminate, $\ih$, equipped with the $\ih$-adic topology. A \textbf{star product} on $\alg$ is an associative $\rea[[\ih]]$-bilinear map $\hstar:\algh \times \algh \to \algh$ having as a unit the constant function $1$, and such that for $A, B \in \alg$,
\begin{align}\label{starproductdefined}
A \hstar B = AB + \sum_{r \geq 1}\ih^{r}B_{r}(A, B),
\end{align}
where the $B_{r}$ are bilinear operators on $\alg$. A star product is called \textbf{symmetric} if $B_{r}(B, A) = (-1)^{r}B_{r}(A, B)$ (which is the same as $A  \hstar B = B \chstar A$ where $\bih = -\ih$), and \textbf{commutative} if $A \hstar B = B \hstar A$. The given associative multiplication on $\alg$ may be extended $\ih$-linearly to all of $\alg[[\ih]]$, and so $\alg[[\ih]]$ has as a topological $\rea[[\ih]]$-module two different associative multiplicative structures which agree to zeroth order in $\ih$. Note that a commutative star product is just a deformation of $\alg$ through commutative, associative algebras

A \textbf{Poisson bracket} on $\alg$ is a Lie bracket, $\pois$, satisfying the condition that for each $A \in \alg$, the adjoint operator, $\{A,\,\}$, is a derivation of $\alg$ as an associative algebra. An $\rea$-linear map is an automorphism of the Poisson algebra $(\alg, \pois)$ if it is simultaneously an associative algebra automorphism and a Lie algebra automorphism. A star product on the Poisson algebra $(\alg, \pois)$ is \textbf{adapted} to $\pois$ provided that $B_{1}(A, B) - B_{1}(B, A) = \{A, B\}$.

If $\alg$ is graded as an associative algebra, a star product on $\alg[[\ih]]$ is called \textbf{graded} if for each  $r \in \nat$, $B_{r}$ is a graded linear operator in the sense that the restriction to $\alg_{k} \times \alg_{l}$ of $B_{r}$ takes values in $\alg_{k+l-r}$. If $\alg$ is a subalgebra of the algebra of smooth functions, $C^{\infty}(N)$, on a smooth manifold, $N$, which is graded as an associative algebra, a star product on $\alg[[\ih]]$ is called \textbf{graded differential} if it is graded and if for each $k, l \in \nat$ the restriction to $\alg_{k}\times \alg_{l}$ of $B_{r}$ is a bidifferential operator of order at most $r$. As A. Astashkevich and R. Brylinski have explained, \cite{Astashkevich-Brylinski}, for equivariant star products the notion of graded differential star products seems to be the appropriate substitute for differential star products, (they call a graded star product what would be here called a symmetric, graded, $\pois$-adapted star product). Note that if the grading of $\alg$ is non-trivial a graded differential star product need not be a differential star product.

A \textbf{homomorphism} of star algebras, $(\alg, \hstar)$ and $(\bar{\alg}, \bhstar)$ will mean a topological $\rea[[\ih]]$-module map, $T:\alg[[\ih]] \to \bar{\alg}[[\ih]]$ such that $T(A \hstar B) = T(A)\bhstar T(B)$. Such a map has the form $T = \sum_{r \geq 0}\ih^{r}T_{r}$ where each $T_{r}$ is the $\ih$-linear extension to $\alg[[\ih]]$ of a $\rea$-linear map $T_{r}:\alg \to \alg$. An isomorphism of star algebras is an invertible homomorphism of star algebras. An isomorphism of star algebras for which $T_{0} = Id_{\alg}$ is called a \textbf{gauge equivalence}. (For background on deformation quantization consult the surveys \cite{Gutt}, \cite{Sternheimer}). A \textbf{geometric isomorphism} of star algebras, $(\alg, \hstar)$ and $(\bar{\alg}, \bhstar)$ will mean a bijective $\rea$-linear map, $T:\alg \to \bar{\alg}$ which when extended $\ih$-linearly to a topological $\rea[[\ih]]$-module map $T: \alg[[\ih]] \to \alg[[\ih]]$ satisfies $T(A \hstar B) = T(A)\bhstar T(B)$. By definition a geometric isomorphism of star algebras is a star algebra isomorphism for which $T_{r} = 0$ if $r \geq 1$. The notion of geometric isomorphism is much more rigid than that of isomorphism (from an algebraic viewpoint this is likely not the `right' notion, but seems adequate for present purposes).

Because the differential of a diffeomorphism acts on the fibers of the cotangent bundle by linear transformations the property of a function on the cotangent bundle being polynomial in each fiber is well-defined. Let $\balg = \polb(T^{\ast}M)$ be the Poisson algebra of smooth, fiberwise-polynomial functions on the cotangent bundle of $M$. An element $a \in \balg$ is a smooth function on $C^{\infty}(T^{\ast}M)$ such that for any $x \in M$ the restriction of $a$ to the fiber $T^{\ast}_{x}M$ is a polynomial; in particular, an element of $\balg$ need not have globally bounded degree on the fibers. Let $\alg = \pol(T^{\ast}M)$ be the subalgebra comprising functions of globally bounded fiberwise degree. The subspaces $\alg^{i} \subset \alg$ comprising functions polynomial in the fibers of degree at most $i$ are well-defined, and each is graded $\alg^{i} = \oplus_{j = 0}^{i}\alg_{j}$, where $\alg_{j}$ comprises functions polynomial in the fibers of degree exactly $j$. The direct sum $\oplus_{i = 0}^{\infty}\alg_{i}$ is isomorphic to $\alg$. The tautological Poisson structure, $\pois$, on $C^{\infty}(T^{\ast}M)$ restricts to both $\balg$ and $\alg$, and the Poisson bracket on $\alg$ is graded in the sense that $\{\alg_{i}, \alg_{j}\} \subset \alg_{i+j-1}$. As an algebra, $\balg$, is identified with the space of sections of the bundle, $S(TM)$, of symmetric algebras of the tangent bundle in such a way that the subspace $\alg_{j}$ is canonically identified with the space of sections of the bundle, $S^{j}(TM)$, of completely symmetric $j$-vectors. The space of sections, $\balg$, of $\oplus_{k = 0}^{\infty}S^{k}(TM)$ and the direct sum, $\oplus_{k = 0}^{\infty}\Gamma(S^{k}(TM)) = \alg$ of the spaces of sections are not the same thing; the latter comprises finite formal sums of tensors whereas the former comprises formal sums of tensors finite at any given point of $M$. When viewed in this way $\alg$ will be called the \textbf{symmetric tensor algebra}, and also written $\sy(M)$. $\balg$ contains $\alg$ as a subalgebra which is proper if $M$ is noncompact. The space of linear differential operators on $C^{\infty}(M)$ having globally bounded order is filtered by the subspaces comprising operators of order bounded by a positive integer, and $\alg$ can be viewed as the associated graded algebra of this filtered algebra. 

A \textbf{projective structure} on $M$ is an equivalence class, $\en$, of torsion-free affine connections such that the unparameterized geodesics of any two representative connections are the same. Two projective structures are equivalent if there is a diffeomorphism of $M$ mapping the geodesics of one structure onto the geodesics of the other structure. The main results of this paper are summarized in the following theorem, the proof of which appears in Section \ref{spsection}.
\begin{theorem}\label{sptheorem}
There is associated to each (curved) projective structure on a smooth $n$-dimensional manifold, $M$, a one-parameter family of graded differential star products, $\hmustar$, on $\alg$ adapted to the tautological Poisson structure, $\pois$, such that in the expansion
\begin{align}\label{speq}
A \hmustar B = AB + \tfrac{\ih}{2}\left(\{A, B\} + (n+1 + 2\mu)\partial C_{1}(A, B)\right) + \sum_{r \geq 2}\ih^{r}B^{\mu}_{r}(A, B),
\end{align} 
\begin{itemize}
\item the restriction to $\alg_{k} \times \alg_{l}$ of each of the bilinear operators $B^{\mu}_{r}$ takes values in $\alg_{k+l-r}$;
\item viewing $\alg_{k}$ as $S^{k}(TM)$, the restriction to $S^{k}(TM) \times S^{l}(TM)$ of each $B_{r}^{\mu}$ is a projectively invariant bilinear differential operator of order at most $r$;
\item each $B_{r}^{\mu}$ is a polynomial in $\mu$ of degree at most $r$; 
\item $\partial C_{1}$ (which is symmetric) is the coboundary of the $1$-cochain for the Hochschild cohomology of $\alg$ with coefficients in $\alg$ given by $C_{1} = \tfrac{1}{n+1 + 2\lie_{\teul}}\circ D$, where $\teul$ is the Euler vector field generating fiber dilations on $T^{\ast}M$; $\tfrac{1}{n+1+2\lie_{\teul}}$ is intepreted as a formal power series in the Lie derivative, $\lie_{\teul}$; and $D$ is the divergence operator on $C^{\infty}(T^{\ast}M)$ determined as in Lemma \ref{divlemma} by a choice of representative $\nabla \in \en$.
\end{itemize} 
These star products satisfy $A \hmustar B = B\,_{-\epsilon}{\star}^{-\mu-n-1}\,\,A$, so that $\hmustar$ is symmetric if and only if $\mu = -(n+1)/2$, in which case the star product is denoted simply $\hstar$. 

For the weight $\mu = -(n+1)/2$, the star products, $\hstar$, associated to two projective structures of the same dimension at least $2$ are geometrically isomorphic if and only if the underlying projective structures are equivalent. In particular, distinct projective structures determine distinct star products. Any two projective structures determine gauge equivalent star products, $\hstar$, so that the gauge equivalence class of the star products so determined depends only on the smooth structure on $M$.

Finally, if $X \in \vect(M)$ is an infinitesimal automorphism of the projective structure, $\en$, then $\lie_{X}$ is an inner derivation of the star algebra $(\alg[[\ih]], \pois, \hstar)$ associated to $\en$.
\end{theorem}
When viewed as operators on $\alg$ the $B_{r}$ are formal pseudo-differential operators rather than differential operators. By Theorem 5.1 of \cite{Duval-ElGradechi-Ovsienko} there is in the projectively flat case a unique symmetric, graded-differential, $\pois$-adapted star product on $\alg$, which must therefore be the star product $\hstar$. Thus $\hstar$ recovers in the flat case the star product described in \cite{Lecomte-Ovsienko}, \cite{Duval-Lecomte-Ovsienko}, \cite{Duval-ElGradechi-Ovsienko}, and \cite{Brylinski}, (note that in these papers a star product by definition is symmetric and satisfies $B_{1} = \tfrac{1}{2}\pois$). That the flat case should have a curved generalization was first suggested by Lecomte and Ovsienko, and is proposed explicitly at the end of \cite{Duval-ElGradechi-Ovsienko}. In the non-flat case it is not even clear how one should formulate a characterization (uniqueness statement) of these star products. On the other hand, even in the flat case the existence proofs given here are different than those of \cite{Lecomte-Ovsienko}. In one dimension it should be possible to relate the associative multiplications defined on the space of modular forms in \cite{Cohen-Manin-Zagier} to some specialization of $\hmustar$; some suggestions in this direction appear in Example \ref{starvectorexample}.

The construction of the projectively invariant star product is based on the construction of what Lecomte - Ovsienko have called a \textbf{projectively invariant quantization map}. Let $\emf^{(i_{1}\dots i_{k})}[\lambda]$ be the bundle on $M$ of symmetric $k$-vectors with $-\tfrac{\lambda}{n+1}$-density coefficients. Denote by $\dop^{k}_{\mu, \lambda} = \Hom(J^{k}(\emf[\mu]), \emf[\mu + \lambda])$ the vector space of differential operators of order at most $k$ mapping sections of the bundle, $\emf[\mu]$, of $-\tfrac{\mu}{n+1}$-densities on $M$, to sections of $\emf[\mu + \lambda]$, and set $\dop_{\mu, \lambda} = \cup_{k = 0}^{\infty}\dop^{k}_{\mu, \lambda}$. Let $\alg_{\lambda, k} = \emf^{(i_{1}\dots i_{k})}[\lambda]$ (both notations will be used) so that $\alg_{\lambda} = \oplus_{k \geq 0}\alg_{\lambda, k}$ is the graded vector space of weighted symmetric polyvectors, which is the associated graded vector space, $\graded \dop_{\mu, \lambda}$ of principal symbols. Call a linear map, $\op: \alg_{\lambda} \to \dop_{\mu, \lambda}$, a \textbf{quantization map}, if for each $k \in \nat$ and for each $a \in \emf^{(i_{1}\dots i_{k})}[\lambda]$ (a $\Bar{\,}$ indicates the space of sections of a vector bundle), $\op(a)$ is a differential operator with principal symbol $a$. On a manifold with projective structure, $\en$, a quantization map will be called \textbf{projectively invariant}, if for each $a \in \emfb^{(i_{1}\dots i_{k})}[\lambda]$, the operator $\op(a)$ is a projectively invariant differential operator. In \cite{Lecomte-Ovsienko}, Lecomte and Ovsienko constructed on any manifold equipped with a flat projective structure a projectively invariant quantization map, and in \cite{Bordemann}, M. Bordemann extended their construction to manifolds equipped with curved projective structures. 

In Theorem \ref{selfadjoint} it is shown that the projectively invariant quantization yields operators which are formally self-adjoint with respect to the canonical pairing between compactly supported $\mu$-densities and $(1-\mu)$-densities. In particular, the representation on operators on half-densities is by formally self-adjoint operators. In \cite{Duval-Lecomte-Ovsienko} this was proved in the flat case by invoking the uniqueness of the flat projectively invariant quantization map proved in \cite{Lecomte-Ovsienko}; here the proof is by integration by parts because the proper formulation of a uniqueness statement is lacking. The star products of Theorem \ref{sptheorem} arise by transporting to $\alg$ via the projectively invariant quantization map the structure of the algebra of differential operators acting on densities. For densities of general weight this yields the associative algebra deformations $\hmustar$ of $\alg$, and when these differential operators act on half-densities the symmetry of the star product, $\hstar$, of Theorem \ref{sptheorem} is a consequence of the self-adjointness of the operators.

The coefficients of a projectively invariant quantization map, $\op(a)$, are linear differential operators in the components of $a$, and so for $f \in \emfb[\mu]$, $\op(a)(f)$ may be viewed as a projectively invariant bilinear differential operator. Lecomte-Ovsienko's and Bordemann's quantization maps are are recovered as special cases of the following general theorem yielding projectively invariant multilinear differential operators. 
\begin{theorem}\label{maintheorem}
Let $(M, \en)$ be a smooth $n$-dimensional manifold with projective structure. Let $\alpha$ be an integer greater than $1$. Let $1 \leq \ga \leq p$; let $k_{\ga}$ be a non-negative integer; let $\lambda_{\ga}$ be a real number not in the set of \textbf{excluded weights}, $\{-n-k_{\ga}, \dots, -n - 2k_{\ga} + 1  \}$; and set $K = \sum_{\ga = 1}^{p}k_{\ga}$ and $\Lambda = \sum_{\ga = 1}^{p}\lambda_{\ga}$. For each integer $1 < \beta \leq K$ there exists a multilinear projectively invariant differential operator of order $K - \beta$:
\begin{align*}
\op_{\be}:\emfb^{(i_{1}\dots i_{k_{1}})}[\lambda_{1}] \times \dots \times \emfb^{(i_{1}\dots i_{k_{p}})}[\lambda_{\al}] \to \emfb^{(i_{1}\dots i_{\beta})}[\Lambda]
\end{align*}
Moreover, if $a \in \emfb^{(i_{1}\dots i_{k})}[\lambda]$ and $f \in \emfb[\mu]$ then the principal symbol of $\op_{0}(a, f)$ viewed as a differential operator in $f$ is $a$. In particular, if $a^{i} \in \emfb^{i}$ is a vector field then $\op_{0}(a, f) = \lie_{a}f$ where $\lie_{a}$ is the usual Lie derivative. 
\end{theorem}
In the multilinear case there is lacking an invariant characterization of these operators, though certainly there should be one. Conceptually the proof of Theorem \ref{maintheorem} is very simple. A projective structure determines on the total space of a $\rea^{\times}$ principal bundle over $M$ a Ricci-flat, torsion-free affine connection, $\hnabla$, the \textbf{Thomas ambient connection}, and for non-excluded weights, $\lambda$, an invariant lift of $a \in \emfb^{(i_{1}\dots i_{k})}[\lambda]$ to a contravariant symmetric $k$-tensor, $\tilde{a}$, of the appropriate homogeneity, and such that $\hnabla \tilde{a}$ is completely trace-free. The operators of Theorem \ref{maintheorem} arise by applying repeatedly the connection $\hnabla$ to the symmetric product of the invariant lifts of the sections of $\emf^{(i_{1}\dots i_{k_{\ga}})}[\lambda_{\ga}]$, and taking as many traces as is possible. This procedure is explicit in the sense that for any particular example it is straightforward to find explicit expressions for the resulting multilinear differential operators. The sense in which Theorem \ref{maintheorem} generalizes the Rankin-Cohen brackets as described in \cite{Cohen-Manin-Zagier} or \cite{Zagier} is explained in Section \ref{onedimsection}. Probably the operators $\op_{\beta}$ generalize the multilinear differential operators on the space of modular forms constructed by R. A. Rankin in \cite{Rankin}, but a precise statement is lacking. As is stated precisely in Corollary \ref{invliftprop}, that the weight $\lambda$ is one of the excluded weights of Theorem \ref{maintheorem} corresponds to the existence of a projectively invariant differential operator on the bundle $\emf^{(i_{1}\dots i_{k})}[\lambda]$. These operators arise as the obstruction to the construction of the invariant lift, and are described in Proposition \ref{droptwoprop}.

The complicated interaction of the invariant lift with the algebra structure of $\alg$ underlies Theorems \ref{sptheorem} and \ref{maintheorem}. Proposition \ref{liftsymprop} shows that the difference between the symmetric product of the invariant lifts of elements of $\sy(M)$ and the invariant lift of the symmetric product of the same elements is expressible explicitly in terms of the operators $\op_{\be}$ of Theorem \ref{maintheorem}. In some sense this statement characterizes these invariant multilinear differential operators. It also leads to the following theorem proved in Section \ref{infsection}.
\begin{theorem}\label{inftheorem}
The limit $\lim_{t \to 0} A \,_{t}{\star}^{c/t}\,\,B$ exists and defines on $\alg$ a commutative, graded differential star product, $\cinfstar$, having the form
\begin{align}\label{speq2}
A \cinfstar B = AB + c\partial C_{1}(A, B) + \sum_{r \geq 2}c^{r}B_{r}^{\infty}(A, B),
\end{align}
where for $1 \leq r$ the restriction to $\alg_{k} \times \alg_{l}$ of $B_{r}^{\infty}$ has the form
\begin{align}\label{binfr}
r!B_{r}^{\infty}(a, b) = \frac{\binom{k+l}{r}}{\binom{n + 2k + 2l - r}{r}}\op_{k+l-r}(a, b).
\end{align}
\end{theorem} 
The star product $\cinfstar$ is new even in the flat case. As is partly proved in Section \ref{infsection}, The real part of the specialization at $c = 2i$ of the limiting star products $\cinfstar$ should recover the commutative multiplication, $\cmzinf$, on the space of modular forms constructed in \cite{Cohen-Manin-Zagier}. A special case of this claim is proved as Example \ref{infmultex}. It is not clear what is the geometric meaning of $\cinfstar$.

All the constructions use in an essential way the ambient viewpoint first developed by T. Y. Thomas for projective and conformal structures, \cite{Thomas}. The work of C. Fefferman on invariants of CR structures gave further impetus to this ambient point of view which has played a prominent role in the construction of local invariants and invariant differential operators on manifolds with curved geometric structures. On a manifold with (curved) conformal structure a construction utilizing instead of the Thomas ambient connection the Fefferman-Graham ambient connection and a conformally invariant lift of tensors (e.g. that used by R. Jenne, \cite{Jenne}), ought to lead to conformally invariant star products. Here much of the underlying work of constructing the right sort of conformally invariant multilinear differential operators on trace-free symmetric polyvectors has been done already by M. Eastwood, \cite{Eastwood}. Also relevant should be the conformally invariant operators on densities constructed by S. Alexakis in \cite{Alexakis}, work which built on the earlier study of the flat case made by Eastwood - C. R. Graham, \cite{Eastwood-Graham}. More generally, it seems plausible that it should be posssible to construct on any manifold equipped with a parabolic geometry an invariant star product. In this context, it would be interesting to know if the projectively invariant star products constructed here can be obtained by transporting one of the standard star products on the symmetric algebra of $\mathfrak{sl}(n+1, \rea)$ using the canonical regular normal Cartan connection associated to the projective structure. In Theorem 3.6 of \cite{Calderbank-Diemer}, D. Calderbank and T. Diemer construct on manifolds with parabolic geometric structure invariant bilinear differential pairings of sections of the vector bundles associated via the Cartan connection to certain Lie algebra cohomologies. In the setting of projective structures it seems possible that their operators are related in some way to those of Theorem \ref{maintheorem}, though this possibility is not explored here. The nature of the construction in \cite{Calderbank-Diemer} enabled its authors to explore in Section 6 of \cite{Calderbank-Diemer} some complicated (curved $A_{\infty}$) algebraic structure related to their operators; an analogous understanding of the algebraic structure of the operators of Theorem \ref{maintheorem} would be desirable.

Theorem \ref{sptheorem} suggests the possibility of an algebraic broadening of the notion of projective structure. The papers, \cite{Connes-Moscovici-modular}, \cite{Connes-Moscovici-rc}, of A. Connes and H. Moscovici, achieve such a goal in the one-dimensional case. Connes and Moscovici define generalized Rankin-Cohen brackets for any associative algebra equipped with an action of the Hopf algebra of transverse geometry in codimension one such that the derivation corresponding to the Schwarzian derivative is inner. Since the Rankin-Cohen brackets arise as a one-dimensional special case of Theorem \ref{maintheorem}, and the Hopf algebras of transverse geometry are defined in higher dimension, it seems plausible that there could be in all dimensions, and in the non-flat case, an algebraic extension of the notion of projective structure along these lines.

\section{Review of Projective Structures}
In this section the basic results concerning projective structures are reviewed, following the approach of T. Y. Thomas, \cite{Thomas}. Some version of much of this material can be found also in various modern sources, for instance \cite{Bailey-Eastwood-Gover} or \cite{Gover}. 

For $a \in \rea$ and $r \in \nat$ write $a_{(r)} = a(a - 1)\dots (a-r+1)$ and $\binom{a}{r} = a_{(r)}/r!$. In tensorial expressions, lowercase Latin indices range over the set $\{1, \dots n\}$ while uppercase Latin indices range over the set $\{\infty, 1, \dots, n\}$. Often the abstract index notation is employed (this should be clear in context), so that expressions with indices have invariant meaning. Complete symmetrization (resp. anti-symmetrization) over a specified set of indices is indicated by enclosing them in parentheses (resp. square brackets). The symmetric product of symmetric polyvectors is denoted $a \odot b$, so that, for example, for vector fields $a^{i}$ and $b^{j}$, $(a\odot b)^{ij} = a^{(i}b^{j)} =\tfrac{1}{2}(a^{i}b^{j} + a^{j}b^{i})$.

For a torsion-free connection, $\nabla$, denote by $R_{ijk}\,^{l}$ its curvature tensor, consistent with the convention $2\nabla_{[i}\nabla_{j]}\alpha_{k} = -R_{ijk}\,^{p}\alpha_{p}$. For an arbitrary torsion-free affine connection, $\nabla$, the tensor $R_{ij} = R_{ipj}\,^{p}$ need not be symmetric. Define the \textbf{projective Schouten tensor}, $P_{ij}$, by $(n-1)P_{ij} = R_{ij} - \frac{2}{n+1}R_{[ij]}$, and the \textbf{projective Weyl tensor}, $B_{ijk}\,^{l} = R_{ijk}\,^{l} + 2\delta_{[i}\,^{l}P_{j]k} - 2P_{[ij]}\delta_{k}\,^{l}$. Tracing the first Bianchi identity gives $2(n+1)P_{[ij]} = 2R_{[ij]} = R_{ijp}\,^{p}$. Torsion-free affine connections have the same unparameterized geodesics if and only if there exists a one-form, $\gamma$, so that their difference tensor has the form $2\gamma_{(i}\delta_{j)}\,^{k}$, in which case $B_{ijk}\,^{l}$ does not depend on the choice of representative $\nabla \in \en$. Tracing the second Bianchi identity gives $\nabla_{[i}P_{jk]} = 0$, so that the two-form $P_{[ij]}$ is closed. The curvature of the $\rea^{\times}$ principal connection induced by $\nabla$ on the principal bundle of frames in $\wedge^{n}(T^{\ast}M)$ is a scalar multiple of $-R_{ijp}\,^{p} = -2(n+1)P_{[ij]}$, which shows that the closed two-form, $-2(n+1)P_{[ij]}$, can be regarded as a scalar multiple of the curvature of the $\rea^{\times}$ principal connection induced by $\nabla$ on this principal bundle. 
If $\nabla\mu = 0$, then $-R_{ijp}\,^{p}\mu_{k_{1}\dots k_{n}} = 2\nabla_{[i}\nabla_{j]}\mu_{k_{1}\dots k_{n}} = 0$, and so by the traced first Bianchi identity, the Ricci tensor of $\nabla$ is symmetric, and so also $P_{[ij]} = 0$. For any $\nabla \in \en$, there is a one-form, $\gamma$, such that $\nabla\mu = (n+1)\gamma\tensor\mu$. If the difference tensor of $\tilde{\nabla} \in \en$ with $\nabla$ is $2\gamma \odot \delta$, then $\tilde{\nabla}\mu = \nabla\mu - (n+1)\gamma\tensor \mu = 0$. This shows that given a projective structure, $\en$, on an orientable manifold, $M$, for each choice of a volume form, $\mu$, on $M$ there exists a unique $\nabla \in \en$ making $\mu$ parallel. Consequently it suffices to work only with representatives $\nabla \in \en$ having symmetric Ricci tensor, and this will be henceforth assumed always.

Let $\mathsf{F}$ denote the $\rea^{\times}$ principal bundle of frames in the canonical bundle, $\wedge^{n}(T^{\ast}M)$, of the smooth $n$-dimensional manifold, $M$. When $n = 2l$ let the $\rea^{\times}$ principal bundle, $\dens$, be the unique $\frac{1}{n+1}$-root of $\mathsf{F}$. When $n = 2l-1$, assume $M$ is orientable with a fixed orientation, so that the group of $\mathsf{F}$ is reduced to $\rea^{>0}$, and let $\dens$ be a choice of $\frac{1}{n+1}$-root of $\mathsf{F}$. Denote by $\emf[\lambda]$ the line bundle associated to $\dens$ by the representation, $r\cdot s = r^{-\lambda}s$, of $\rea^{\times}$ on $\rea$, so that $\emf[-1]$ is a $\frac{1}{n+1}$-root of $\Lambda^{n}(T^{\ast}M)$. A section of $\emf[\lambda]$ is a $-\lambda/(n+1)$ density. Sections, $u$, of $\emf[\lambda]$ are in canonical bijection with functions $\tilde{u}:\dens \to \rea$ homogeneous of degree $\lambda$. The model for $\dens$ is the the defining bundle $\standrepz  = \standrep - \{0\}\to \proj(\standrep)$, where $(\standrep, \svol)$ is an $n+1$ dimensional real vector space with volume form $\svol$. This defining bundle is the bundle of frames in the tautological line bundle over $\proj(\standrep)$. Write $\emf^{i_{1}\dots i_{k}}_{j_{1}\dots j_{l}}[\lambda] = \T^{k}(TM) \tensor \T^{l}(T^{\ast}M) \tensor \emf[\lambda]$. The use of $(\,)$ and $[\,]$ indicates symmetries, so that, for instance, $\emf^{(ij)k}[\lambda]$ indicates weighted $3$-vectors symmetric in the first two indices.

Fix a choice of $\rho:\dens \to M$. A tautological $n$-form, $\alpha$, is defined on $\dens$, by 
\begin{align*}
&\alpha_{s}(X_{1}, \dots, X_{n}) = s^{n+1}_{\rho(s)}(\rho_{\ast}(X_{1}), \dots, \rho_{\ast}(X_{n}))& & \text{for}\, X_{1}, \dots, X_{n} \in T_{s}\dens.&
\end{align*}
A canonical volume form, $\svol$, is defined by $\svol = d\alpha$. A canonical {\bf Euler vector field}, $\eul$, the infinitesimal generator of the fiber dilations on $\dens$, is defined by $i(\eul)\svol = (n+1)\alpha$.

A torsion-free affine connection, $\hat{\nabla}$, on $\dens$, and a choice of section, $s$, determine on $M$ a torsion-free affine connection, $\nabla$, defined by $\nabla_{X}Y = \rho_{\ast}(\Hat{\nabla}_{\hat{X}}\hat{Y})$. If the condition $\hat{\nabla}\eul = \delta$ is imposed, it is easily checked that the connection, $\tilde{\nabla}$, determined by $\tilde{s} = fs$ has the same unparameterized geodesics as does $\nabla$.
This associates a projective structure on $M$ to each torsion-free affine connection, $\hat{\nabla}$, on $\dens$, satisfying $\hat{\nabla}\eul = \delta$. It may be checked that if $\hat{\nabla}\svol = 0$ then $\nabla$ is the unique representative of $\en$ making $\mu= s^{\ast}(\alpha)$ parallel. Two torsion-free connections, $\hat{\nabla}$ and $\hat{\nabla}^{\prime}$, on $\dens$, satisying $\hat{\nabla}\eul = \delta = \hat{\nabla}^{\prime}\eul$ and making parallel $\svol$ determine on $M$ the same projective structure. The only freedom is in the vertical part of $\hat{\nabla}_{\hat{X}}\hat{Y}$, and requiring that $\hat{\nabla}$ be Ricci flat eliminates this freedom.

\begin{theorem}[T. Y. Thomas, \cite{Thomas}]\label{thomastheorem}
Let $M$ be a smooth manifold equipped with a choice of $\rho:\dens\to M$. There is a functor associating to each projective structure, $\en$, on $M$, a unique torsion-free affine connection, $\hnabla$, (the Thomas ambient connection) on $\dens$ satisfying the following conditions:
\setcounter{condition}{0}
\begin{condition}\label{tcon1}
$\hnabla\eul$ is the fundamental $\binom{1}{1}$-tensor on $\dens$. 
\end{condition}
\begin{condition}\label{tcon2}
 $\hnabla\svol = 0$.
\end{condition}
\begin{condition}\label{tcon3}
$\hnabla$ is Ricci flat.
\end{condition}
\begin{condition}\label{tcon4}
The projections into $M$ of the unparametrized geodesics of $\hnabla$ transverse to the vertical are the geodesics of the given projection structure.
\end{condition}
\end{theorem}
 The choice of a local section, $s$, induces on $\dens$, a unique $\rea^{\times}$ principal connection, $\phi$, such that $s$ is a parallel section and $\phi(\eul) = 1$.  The connection, $\phi$, determines a horizontal lift, $\hat{X}$, of each vector field, $X$, on $M$. Condition \ref{tcon4} has the following analytic reformulation.
\begin{condition}\label{tcon4b}
For each (local) section, $s$, of $\rho:\dens\to M$ the connection, $\nabla$, defined on $M$ by $\nabla_{X}Y = \rho_{\ast}(\hnabla_{\hat{X}}\hat{Y})$, is the unique representative of $\en$ making the volume form $s^{\ast}(\alpha)= s^{n+1}$ parallel.
\end{condition}
\noindent
The theorem is proved first for a neighborhood on $M$ for which $\dens$ is trivial. The uniqueness statement in the theorem guarantees that the resulting connections patch together globally on $\form$. To sketch the proof, suppose there is a global section, $s:M \to \dens$. Let $\nabla \in \en$ be the unique representative of the given projective structure $\nabla \in \en$ making $\mu = s^{\ast}(\alpha)$ parallel. For any symmetric tensor, $P_{ij}$, conditions \eqref{tcon1}-\eqref{tcon3} of Theorem \ref{thomastheorem} are satisfied by the connection, $\hnabla$, defined by requiring it to be torsion-free and to satisfy
\begin{align*}
&\hnabla_{V}\eul = V,& &\hnabla_{\hat{X}}\hat{Y} = \widehat{\nabla_{X}Y} + P(X, Y)\eul.&
\end{align*}
Straightforward computation of the curvature of $\hnabla$ shows that requiring $\hnabla$ to be Ricci-flat determines $P_{ij}$ uniquely as $P_{ij} = \frac{1}{n-1}R_{ij}$, where $R_{ij} = R_{ipj}\,^{p}$ is the Ricci tensor of $\nabla$. 

Because the curvature tensor of $\hnabla$ is horizontal, its components may be regarded as tensors on $M$. The possibly non-vanishing components are $B_{ijk}\,^{l}$ and $C_{ijk} = 2\nabla_{[i}P_{j]k}$. The Bianchi identities imply $\nabla_{[i}B_{jk]l}\,^{p} = -\delta_{[i}\,^{p}C_{jk]l}$ and $(2-n)C_{ijk} = \nabla_{p}B_{ijk}\,^{p}$, so that if $n > 2$ the vanishing of $B_{ijk}\,^{l}$ implies the vanishing of $C_{ijk}$ and consequently of the curvature of $\hnabla$. If $n = 2$, $B_{ijk}\,^{l} = 0$ automatically, and $C_{ijk} = 0$ implies $\hnabla$ is flat. The Thomas ambient connection contains exactly the same information as does the canonical regular, normal Cartan connection inducing the projective structure, and each may be recovered from the other straightforwardly.

 By considering the transformation law for $P_{ij}$ under a change of scale, and considering formal power series, it is not hard to show (see \cite{Gover}) that there exists at each $x \in M$ a \textbf{projective normal scale}, namely a projective scale, $s$, for which at $x$ there holds $\nabla_{(i_{1}\dots i_{k-2}}P_{i_{k-1}i_{k})} = 0$ for all $k \geq 2$.

The affine space, $\aff(M)$, of projective structures on $M$ is modeled on the space of trace-free sections of $\affdiff$. The difference tensor of two projective structures is defined by the observation that the difference tensor of the unique representatives of the projective structures making parallel a given projective scale does not depend on the choice of scale. The action of the group, $\diff(M)$, of diffeomorphisms of $M$ on $\aff(M)$ is defined by pullback of differential operators as $\phi \cdot \en = [\phi^{\ast} \circ \nabla \circ (\phi^{\ast})^{-1}]$. Differentiating this action along a one-parameter subgroup yields the action on $\aff(M)$ of the Lie algebra, $\vect(M)$, of vector fields on $M$, as the Lie derivative, $\lie_{X}\en$, of a projective structure along a vector field, $X$. In terms of a representative connection $\nabla \in \en$, $\lie_{X}\nabla$ equals the trace-free part of $\nabla_{i}\nabla_{j}X^{k} - X^{p}R_{pij}\,^{k}$, which is easily checked to be a projectively invariant differential operator. The difference $\phi \cdot \en - \en$ is a non-trivial cocycle of $\diff(M)$, sometimes called the \textbf{Schwarzian derivative} (determined by $\en$) of $\phi$. The Schwarzian cocycles associated to different choices of projective structure are cohomologous, and so there is a canonical Schwarzian cohomology class determined solely by the smooth structure on $M$. This seems related to the conclusion of Theorem \ref{gaugetheorem} that the projectively invariant star products determined by different projective structures are gauge equivalent.

Each (local) section, $s$, of $\form$ determines a trivialization of $\form$; a horizontal lift of vector fields, $X \to \hat{X}$ from $M$ to $\form$; a volume form $\mu = s^{n+1}$; and a unique representative $\nabla \in \en$ such that $\nabla \mu = 0$. Each $\phi \in \diff(M)$ has a canonical lift to a principal bundle automorphism $\tilde{\phi} \in \aut(\form)$, and this induces a canonical lift of vector fields associating to each $X \in \vect(M)$ the homogeneity $0$ vector field on $\form$ expressible as
\begin{align}\label{invariantliftvector}
\tilde{X} = \hat{X} - \tfrac{1}{n+1}\diver_{\mu}(X)\eul = \tilde{X} = \hat{X} - \tfrac{1}{n+1}\tr(\nabla X)\eul.
\end{align}
$\tilde{X}$ is determined uniquely by the requirements that $\rho_{\ast}(\tilde{X})= X$, $\lie_{\eul}\tilde{X} = 0$, and $\lie_{\tilde{X}}\Psi = 0$ (Equivalently $\tr \hnabla \tilde{X} = 0$). While the invariant lift is a linear map, it is not a $C^{\infty}$-module homomorphism; precisely, $\widetilde{fX} = f\tilde{X} - \tfrac{1}{n+1}df(X)\eul$. As a consequence, the invariant lift of vector fields does not induce directly an invariant lift of polyvectors. However expressing the identity $\lie_{A}\Psi = 0$ in terms of the ambient connection as $\tr \hnabla A = 0$ suggests that a projective structure should determine an invariant lift of completely symmetric polyvectors as follows. Call $\tilde{a} \in S^{k}(T\form)$ a lift of $a \in \sklb$ if for any $\theta^{1}, \dots, \theta^{k} \in \Gamma(T^{\ast}M)$, $\tilde{a}(\rho^{\ast}(\theta^{1}), \dots, \rho^{\ast}(\theta^{k}))$ is the homogeneity $\lambda$ function on $\form$ corresponding to the density $a(\theta^{1}, \dots, \theta^{k})$. The invariant lift, $\tilde{a}$, of $a \in \emf^{(i_{1}\dots i_{k})}$ is determined uniquely by the requirements $\lie_{\eul}\tilde{a} = 0$ and $\tr \hnabla \tilde{a} = 0$. More generally, it is possible to lift weighted symmetric polyvectors, except for a discrete set of excluded weights which correspond to the existence of projectively invariant differential operators obtained by taking appropriate combinations of the complete traces of some number of covariant derivatives of a section of $\skl$. The existence of the lift described in Proposition \ref{invliftprop} was proved already by M. Bordemann in \cite{Bordemann}.

\begin{proposition}\label{invliftprop}
Let $a \in \sklb$. If $\lambda \notin \{-n-k-m: 0 \leq m \leq k-1\}$
%\{-n - k, -n - k - 1, \dots, -n -2k +1\}$
 there is a projectively invariant lift, $\tilde{a} \in \Gamma(S^{k}(T\form))$, characterized uniquely by the requirements $\lie_{\eul}\tilde{a} = \lambda \tilde{a}$ and $\tr \hnabla \tilde{a} = 0$.

For each excluded weight, $\lambda = -n-k-m$, there exists a projectively invariant differential operator, $\kop:\emf^{(i_{1}\dots i_{k})}[\lambda] \to \emf^{(i_{1}\dots i_{m})}[\lambda]$, so that for any $\nabla \in \en$, $\kop$ has the form
\begin{align}\label{invinter}
\kop(a)^{i_{1} \dots i_{m}} = \nabla_{j_{1}}\dots \nabla_{j_{k-m}}a^{i_{1}\dots i_{m}j_{1}\dots j_{k-m}} + \sum_{s = 2}^{k-m}C^{s}_{i_{1}\dots i_{s}}\nabla_{j_{1}}\dots \nabla_{j_{k-m-s}}a^{i_{1}\dots i_{s}j_{1}\dots j_{k-m-s}},
\end{align}
where $C^{s}$ is a polynomial in $P_{ij}$ and the completely symmetrized covariant derivatives of $P_{ij}$ of order at most $s - 2$.
\end{proposition}
\begin{proof}
In this proof, as opposed to in the rest of the paper, a local frame will be chosen and indices will refer to the chosen frame and dual coframe. Let $\nabla \in \en$ be the representative associated to a choice of local section $s:M \to \form$. Let $E_{i}$ and $\theta^{j}$ be a local frame and dual coframe on the base, and take a coframe $\phi^{I}$ on $\form$ where $\phi^{i} = \rho^{\ast}(\theta^{i})$ and $\phi^{\infty}$ is the principal connection determined by the scale, $s$ (so that on $\form$ a dual frame is given by $F_{i} = \hat{E}_{i}$ and $F_{\infty} = \eul$). The first two requirements mean that if $\tilde{a}$ is to exist there exist $a_{m} = a_{m}^{i_{1}\dots i_{m}}E_{i_{1}}\odot \dots \odot E_{i_{m}} \in \emfb^{(i_{1}\dots i_{m})}[\lambda]$ such that $\tilde{a}$ can be written in the form
\begin{align}\label{undeterminedlift}
&\tilde{a} = \sum_{m = 0}^{k}\tilde{a}_{m}^{i_{1}\dots i_{m}}\hat{E}_{i_{1}}\odot \dots \odot \hat{E}_{i_{m}}\odot \eul^{k-m},& &\eul^{k-m} = \eul \odot \dots \odot \eul,&
\end{align}
where $\tilde{a}_{m}^{i_{1}\dots i_{m}}$ is the homogeneity $\lambda$ function on $\form$ corresponding to the components $a_{m}^{i_{1}\dots i_{m}}$, which are viewed as densities on $M$. The condition $\tr\hnabla \tilde{a} = 0$ gives for each $1 \leq m \leq k-1$ a linear equation relating $a_{m-1}$, $a_{m}$, and $a_{m + 1}$, and for $m = k$ a linear equation relating $a_{m-1}$ and $a_{m}$. If $-\lambda \notin \{n + k, n + k + 1, \dots, n +2k -1\}$ these equations may be solved for $a_{m-1}$. Explicitly one obtains the identities:
\begin{align}
&\label{induct1}(\lambda + n + 2k-1)a_{k-1}^{i_{1}\dots i_{k-1}} = -k\nabla_{p}a_{k}^{i_{1}\dots i_{k-1} p}, \qquad \text{and, for\,\,} 1 \leq m \leq k -1,\\
&\label{induct2}(\lambda + n + m + k - 1)a_{m - 1}^{i_{1}\dots i_{m - 1}} = -\tfrac{m}{k-m + 1}(\nabla_{p}a_{m}^{i_{1}\dots i_{m - 1}p} + (m + 1)a_{m + 1}^{i_{1}\dots i_{m -1}pq}P_{pq}).& 
\end{align}
To obtain \eqref{induct1} and \eqref{induct2} proceed as follows. Write $\tilde{a}$ as in \eqref{undeterminedlift}. Then $\binom{k}{m}\tilde{a}^{i_{1}\dots i_{m}\infty \dots \infty} = \tilde{a}_{m}^{i_{1}\dots i_{m}}$. As $\hnabla_{\eul}\tilde{a} = (\lambda + k)\tilde{a}$, there follows
\begin{align}\label{intermed0}
\binom{k}{m}\hnabla_{\infty}\tilde{a}^{i_{1}\dots i_{m}\infty \dots \infty} =(\lambda + k)\tilde{a}_{m}^{i_{1}\dots i_{m}},
\end{align}
where the expression on the left hand side refers to components of $\hnabla_{\eul}\tilde{a}$ with respect to the chosen frame and dual coframe. Likewise careful handling of the symmetrization operator yields
\begin{align*}
&\hnabla_{\hat{E}_{i}}\tilde{a} = 
\sum_{m = 0}^{k}\hnabla_{\hat{E}_{i}}(\tilde{a}_{m}^{i_{1}\dots i_{m}})\hat{E}_{i_{1}}\odot \dots \odot \hat{E}_{i_{m}}\odot \eul^{k-m}
 + \sum_{m = 1}^{k}m\tilde{a}_{m}\,^{i_{1}\dots i_{m}}(\widehat{\nabla_{E_{i}}E_{i_{1}}})\odot \dots \odot \hat{E}_{i_{m}}\odot \eul^{k-m}\\& 
+ \sum_{m = 1}^{k}m\tilde{a}_{m}\,^{i_{1}\dots i_{m-1}p}P_{ip}\hat{E}_{i_{1}}\odot \dots \odot \hat{E}_{i_{m-1}}\odot \eul^{k-m+1} \\&
+ \sum_{m =0}^{k-1}(k-m)\tilde{a}_{m}^{(i_{1}\dots i_{m}}\delta_{i}\,^{i_{m +1})}\hat{E}_{i_{1}}\odot \dots \odot \hat{E}_{i_{m + 1}}\odot \eul^{k-m-1}.
\end{align*}
It will make formulas easier to read if notation is abused in that there is dropped in expressions such as $\tilde{a}_{m}^{i_{1}\dots i_{m}}$ the $\tilde{\,\,}$, so that there is written, for instance, $\binom{k}{m}\tilde{a}^{i_{1}\dots i_{m}\infty \dots \infty} = a_{m}^{i_{1}\dots i_{m}}$. The preceeding computation of $\hnabla_{\hat{E}_{i}}\tilde{a}$ implies the identities
\begin{align}
&\label{intermed1}\hnabla_{i}\tilde{a}^{i_{1}\dots i_{k}} = \nabla_{i}a_{k}^{i_{1}\dots i_{k}} + a_{k-1}^{(i_{1}\dots i_{k-1}}\delta_{i}^{i_{k})}, \qquad \text{and, for\,\,} 1 \leq m \leq k - 1,\\
&\label{intermed2}\binom{k}{m}\hnabla_{i}\tilde{a}^{i_{1}\dots i_{m}\infty \dots \infty} = \nabla_{i}a_{m}^{i_{1}\dots i_{m}} + (k - m + 1)a_{m - 1}^{(i_{1}\dots i_{m - 1}}\delta_{i}^{i_{m})} + (m + 1)a_{m + 1}^{i_{1}\dots i_{m}p}P_{ip}.
\end{align}
Taking a trace and manipulating coefficients gives \eqref{induct1} and \eqref{induct2}. This completes the proof of the existence of the invariant lift for non-excluded weights.

If $\lambda = -n - k - m$ with $m \in \{0, \dots, k-1\}$ then one can look for $\tilde{a}$ of the form \eqref{undeterminedlift} and satisfying the condition $\lie_{\eul}\tilde{a} = \lambda \tilde{a}$. Imposing for $m < \gamma < k$ the conditions $\hnabla_{Q}\tilde{a}^{i_{1}\dots i_{\gamma}Q\infty\dots \infty} = 0$, gives equations of the form \eqref{induct1} and \eqref{induct2} which determine uniquely the components $a_{\gamma}^{i_{1}\dots i_{\gamma}}$ for $\gamma > m$. Moreover, the components $\hnabla_{Q}\tilde{a}^{i_{1}\dots i_{m}Q\infty\dots \infty}$ are determined uniquely as is seen easily using \eqref{intermed1} and \eqref{intermed2}; solving inductively for the components $a_{m}^{i_{1}\dots i_{m}}$ in terms of $a^{i_{1}\dots i_{k}}$ shows that they are expressible as an invariant differential operator of order $k-m$ applied to $a^{i_{1}\dots i_{k}}$; having top order piece a constant multiple of the top order piece of $\nabla_{j_{1}}\dots \nabla_{j_{k-m}}a^{i_{1}\dots i_{m}j_{1}\dots j_{k - m}}$; and having terms of order $0 \leq s \leq k-m-2$ as in \eqref{invinter} which are expressible as polynomials in $P_{ij}$ and its covariant derivatives of order at most $s$. Since any appearance of $P_{ij}$ and its covariant derivatives is contracted completely with a completely symmetric tensor built from $a$ and its covariant derivatives, any covariant derivative of $P_{ij}$ may be replaced by its complete symmetrization.
\end{proof}

\begin{remark}
In the case of excluded weights the conditions imposed in the proof of Proposition \ref{invliftprop} do not determine the components $a_{s}^{i_{1}\dots i_{s}}$ for $s \leq m$, and these components may be chosen arbitrarily if it is desired to construct some lift of $a$. Such a lift will be invariant modulo $\eul^{k-m}$.
\end{remark}

\begin{remark}
In general there should exist an invariant lift of weighted contravariant tensors having symmetries described by a given Young diagram except for a set of excluded weights that can be read off from the Young diagram. Excluded weights should correspond to the existence of projectively invariant differential operators. For example, except for the excluded weight $\lambda = -n - 2k +1$, $a^{i_{1}\dots i_{k}} \in \emfb^{[i_{1}\dots i_{k}]}[\lambda]$ admits the invariant lift 
\begin{align*}
\tilde{a} = \widetilde{a^{i_{1}\dots i_{k}}}\hat{E}_{i_{1}}\wedge \dots \wedge \hat{E}_{i_{k}} - \tfrac{k}{\lambda + n + 2k-1}\widetilde{\nabla_{p}a^{i_{1}\dots i_{k-1}p}}\hat{E}_{i_{1}}\wedge \dots \hat{E}_{i_{k-1}}\wedge \eul,
\end{align*}
where here the wedge product means complete antisymmetrization of the tensor product. The operator $\kop:\emfb^{[i_{1}\dots i_{k}]}[\lambda] \to \emfb^{[i_{1}\dots i_{k-1}]}[\lambda]$ defined by $\kop(a) = \nabla_{p}a^{i_{1}\dots i_{k-1}p}$ is projectively invariant if and only if $\lambda = -n-2k+1$. 
\end{remark}

\begin{remark}
Something like the invariant lift has appeared in previous papers treating conformal and projective invariants and invariant differential operators. In the note \cite{Thomas-conformalnote1}, Thomas defined a conformally invariant lift of weighted covariant skew-symmetric two-tensors (he called the lift `completing the tensor'). In his thesis, \cite{Jenne}, Jenne used the conformally invariant lift of symmetric tensors to construct conformally invariant multilinear differential operators on weighted tensors. In eq. $(33)$ of \cite{Gover}, R. Gover uses a projectively invariant lift of an unweighted tensor of curvature tensor type. Similar constructions have been used, usually implicitly, in various works exploiting the ambient point of view, e.g. by Eastwood, C. Fefferman, Gover, Graham, K. Hirachi, etc. . . 
\end{remark}

\begin{example}
When $k = 1$, \eqref{induct1} gives $(\lambda + n + 1)a_{0} = - \nabla_{p}a_{1}^{p}$. When $\lambda = 0$, this recovers \eqref{invariantliftvector}. The role of the excluded weight $-n-1$ is indicated by the observation that the operator $\kop:\emfb^{i}[\lambda] \to \emfb[\lambda]$ given by $\kop(a^{i}) = \nabla_{p}a^{p}$ is projectively invariant if and only if $\lambda = -n-1$. 
\end{example}
\begin{example}
When $k = 2$, \eqref{induct1} and \eqref{induct2} give
\begin{align*}
  &(\lambda + n + 3)a_{1}^{i} = -2\nabla_{p}a_{2}^{ip}& &-2(\lambda + n + 2)a_{0} = \nabla_{p}a_{1}^{p} + 2a_{2}^{pq}P_{pq},
\end{align*}
so that if $\lambda \notin \{-n-2, -n-3\}$ the invariant lift of $a^{ij} \in \emfb^{(ij)}[\lambda]$ is defined by $a_{1}^{i} = -\tfrac{2}{\lambda + n + 3}\nabla_{p}a^{ip}$ and $a_{0} = \tfrac{1}{(\lambda + n + 3)_{(2)}}\nabla_{p}\nabla_{q}a^{pq} - \tfrac{1}{(\lambda + n + 2)}a^{pq}P_{pq}$. The excluded weights have the following significance. The differential operator $a^{ij} \to \nabla_{p}a^{ip}$ defined on $\emf^{(ij)}[\lambda]$ is invariant if and only if $\lambda = -n - 3$, and the differential operator $a^{ij} \to \nabla_{p}\nabla_{q}a^{pq} - a^{pq}P_{pq}$ defined on $\emf^{(ij)}[\lambda]$ is invariant if and only if $\lambda = -n - 2$.
\end{example}
\begin{example}
When $k = 3$, \eqref{induct1} and \eqref{induct2} give
\begin{align*}
&(\lambda + n + 5)a_{2}^{ij} = -3\nabla_{p}a^{ijp},&
 &(\lambda + n + 4)a_{1}^{i} = -\nabla_{p}a_{2}^{ip} - 3a^{ipq}P_{pq},&\\
&(\lambda + n + 3)a_{0} = -\tfrac{1}{3}(\nabla_{p}a_{1}^{p} + 2a_{2}^{pq}P_{pq}),
\end{align*}
so that if $\lambda \notin \{-n-3, -n-4, -n -5\}$ the invariant lift of $a^{ijk} \in \emfb^{(ijk)}[\lambda]$ is defined by 
\begin{align*}
&a_{2}^{ij} = -\tfrac{3}{\lambda + n + 5}\nabla_{p}a^{ijp},&\\
&a_{1}^{i} = \tfrac{3}{(\lambda + n + 5)_{(2)}}\nabla_{p}\nabla_{q}a^{ipq} - \tfrac{3}{(\lambda + n + 4)}a^{ipq}P_{pq},\\
&a_{0} = \tfrac{-1}{(\lambda + n +5)_{(3)}}\nabla_{p}\nabla_{q}\nabla_{r}a^{pqr} + \tfrac{1}{(\lambda + n +4)_{(2)}}a^{pqr}\nabla_{p}P_{qr} + \tfrac{1}{(\lambda + n + 5)}(3 + \tfrac{1}{(\lambda + n + 4)})P_{qr}\nabla_{p}a^{pqr}
\end{align*}
The following differential operators corresponding to the excluded weights are invariant:
\begin{align*}
&a^{ijk} \to \nabla_{p}a^{ijp}& &\text{for\,\,} a^{ijk} \in \emfb^{(ijk)}[-n-5],\\
&a^{ijk} \to \nabla_{p}\nabla_{q}a^{ipq} - P_{pq}a^{ipq}& &\text{for\,\,} a^{ijk} \in \emfb^{(ijk)}[-n-4],\\
&a^{ijk} \to \nabla_{p}\nabla_{q}\nabla_{r}a^{pqr} -4P_{pq}\nabla_{r}a^{pqr} - 2a^{pqr}\nabla_{p}P_{qr}& &\text{for\,\,} a^{ijk} \in \emfb^{(ijk)}[-n-3].
\end{align*}
\end{example}

.

\begin{corollary}[Corollary of Proposition \ref{invliftprop}]
If $\en$ admits a Ricci-flat representative, $\nabla$, then for $a \in \emfb^{(i_{1}\dots i_{k})}[\lambda]$ (where $\lambda$ is not an excluded weight)
\begin{align}\label{ricflatlift}
a_{k - m}^{i_{1}\dots i_{k-m}} = \tfrac{(-1)^{m}\binom{k}{m}}{(\lambda + n + 2k-1)_{(m)}}\nabla_{j_{1}}\dots \nabla_{j_{m}}a^{i_{1}\dots i_{k-m}j_{1}\dots j_{m}}.
\end{align}
For a general $\en$, if $\nabla \in \en$ is the representative associated to a projective normal scale at $x \in M$ then at the point $x$ there holds \eqref{ricflatlift}.
\end{corollary}

\begin{proof}
When $\en$ admits a Ricci-flat representative it is straightforward to solve \eqref{induct1} and \eqref{induct2} to obtain explicit formulas for the $a_{m}$. This appears already in M. Bordemann's \cite{Bordemann}. For general $\en$ solving \eqref{induct1} and \eqref{induct2} yields an expression for $a_{m}$ which differs from \eqref{ricflatlift} only by a sum of terms expressible of the form a contraction of the covariant derivative of $a$ of order at most $0 \leq r \leq k-m-2$ with a covariant derivative of $P_{ij}$ of order $r$. Since the coefficient of any appearance of covariant derivatives of $P_{ij}$ is completely symmetric, such a term vanishes in a projective normal scale at $x$. 
\end{proof}

\section{Projectively Invariant Multilinear Differential Operators}
\begin{proof}[Proof of Theorem \ref{maintheorem}]
Define $\op_{\beta}$ by 
\begin{align}\label{opdefined}
\widetilde{\op_{\beta}(a_{1}, \dots, a_{\p})}^{i_{1}\dots i_{\be}} = \hnabla_{P_{1}}\dots \hnabla_{P_{K-\beta}}(\tilde{a}_{1} \odot \dots \odot \tilde{a}_{p})^{(i_{1}\dots i_{\be} P_{1}\dots P_{K - \beta})}.
\end{align}
The special case $\op_{0}:\emf^{(i_{1}\dots i_{k})}[\lambda] \tensor \emf[\mu] \to \emf[\lambda + \mu]$ has the form
\begin{align}\label{pinv}
  \widetilde{\op_{0}(a, f)} = \hnabla_{Q_{1}}\dots \hnabla_{Q_{k}}(\tilde{f}\tilde{a})^{Q_{1}\dots Q_{k}} = \tilde{a}^{Q_{1}\dots Q_{k}}\hnabla_{Q_{1}}\dots \hnabla_{Q_{k}}\tilde{f},
\end{align}
where the second equality follows from $\tr \hnabla \tilde{a} = 0$. From this and the explicit expressions for the components $\tilde{a}^{Q_{1}\dots Q_{n}}$ it is apparent that the principal symbol of $\op_{0}(a, f)$, viewed as a differential operator on $\emf[\mu]$, is $a^{i_{1}\dots i_{k}}$. If $a^{i} \in \emfb^{i}[\lambda]$, then
\begin{align*}
\op_{0}(a, f) = a^{p}\nabla_{p}f - \tfrac{\mu}{(n+1 +\lambda)}f\nabla_{p}a^{p}.
\end{align*}
In particular, if $\lambda = 0$, this shows $\op_{0}(a, f) = \lie_{a}f$.
\end{proof}

\begin{remark}
A construction using complete contractions as in \eqref{opdefined} has been used repeatedly in the construction of conformal and projective invariants. The idea goes back to H. Weyl, \cite{Weyl}. Examples of other papers utilizing such a construction are \cite{Fefferman-Graham}, \cite{Gover}, \cite{Jenne}, \cite{Bailey-Eastwood-Graham}, \cite{Alexakis}, and \cite{Fefferman-Hirachi} (this list is by no means complete).
\end{remark}

The definition \eqref{pinv} is constructive in the sense that it is in principle straightforward to compute the complicated explicit expression for $\op_{\beta}$. To illustrate this some examples are given now. It will be useful to note that, for $f \in \emfb[\mu]$, 
\begin{align}\label{iinfmix}
\eul^{J_{1}}\dots \eul^{J_{l}}\hnabla_{J_{1}}\dots \hnabla_{J_{l}}\hnabla_{I_{1}}\dots \hnabla_{I_{k}}\tilde{f} =  (\mu - k)_{(l)}\hnabla_{I_{1}}\dots \hnabla_{I_{k}}f.
\end{align}

\begin{example}
For $s = 1, 2$ let $a_{s} \in \emfb^{(i_{1}\dots i_{k_{s}})}[\lambda_{s}]$ and write $K = k_{1} + k_{2}$. Then
\begin{align}\label{opkminusone}
&\op_{K-1}(a_{1}, a_{2})^{i_{1}\dots i_{K-1}} = \tfrac{k_{1}}{K}a_{1}^{p(i_{1}\dots i_{k_{1}-1}}\nabla_{p}a_{2}^{i_{k_{1}}\dots i_{K-1})} +  \tfrac{k_{2}}{K}a_{2}^{p(i_{1}\dots i_{k_{2}-1}}\nabla_{p}a_{1}^{i_{k_{2}}\dots i_{K-1})} \\
&\notag - \tfrac{k_{2}(\lambda_{1} + 2k_{1})}{K(\lambda_{2} + n + 2k_{2} - 1)}a_{1}^{(i_{1}\dots i_{k_{1}}}\nabla_{p}a_{2}^{i_{k_{1}+1} \dots i_{K-1})p}   - \tfrac{k_{1}(\lambda_{2} + 2k_{2})}{K(\lambda_{1} + n + 2k_{1} - 1)}a_{2}^{(i_{1}\dots i_{k_{2}}}\nabla_{p}a_{1}^{i_{k_{2}+1} \dots i_{K-1})p}.
\end{align}
To show \eqref{opkminusone} observe that, by definition,
\begin{align*}
&\widetilde{\op_{K-1}(a_{1}, a_{2})}^{i_{1}\dots i_{K-1}} = \hnabla_{I}\tilde{a}_{1}^{(i_{1}\dots i_{k_{1}}}\tilde{a}_{2}^{i_{k_{1}+1}\dots i_{K-1} I)}\\
& = \tfrac{k_{1}}{K}\tilde{a}_{1}^{I(i_{1}\dots i_{k_{1}-1}}\hnabla_{I}\tilde{a}_{2}^{i_{k_{1}}\dots i_{K-1})} +\tfrac{k_{2}}{K}\tilde{a}_{1}^{I(i_{1}\dots i_{k_{2}-1}}\hnabla_{I}\tilde{a}_{2}^{i_{k_{2}}\dots i_{K-1})}. 
\end{align*}
This is straightforward to evaluate explicitly using \eqref{induct1}, \eqref{intermed0}, and \eqref{intermed1}, and the result is \eqref{opkminusone}. For $a \in \emfb^{(i_{1}\dots i_{k_{1}})}[\lambda_{1}]$ and $b \in  \emfb^{(i_{1}\dots i_{k_{2}})}[\lambda_{2}]$, examples of \eqref{opkminusone} are 
\begin{align*}
&\op_{K-1}(a, b)^{i_{1}\dots i_{K-1}}& &k_{1}& &k_{2}\\
\hline\\
&a^{i_{1}\dots i_{s}p}\nabla_{p}b - \tfrac{\lambda_{2}}{(n+2s - 1+\lambda_{1})}b\nabla_{p}a^{pi_{1}\dots i_{s}}  & &s& &0\\
&\tfrac{1}{2}a^{p}\nabla_{p}b^{i} + \tfrac{1}{2}b^{p}\nabla_{p}a^{i} - \tfrac{(\lambda_{1} + 2)}{2(n+2+\lambda_{2})}a^{i}\nabla_{p}b^{p} -  \tfrac{(\lambda_{2} + 2)}{2(n+2+\lambda_{1})}b^{i}\nabla_{p}a^{p}  & &1& &1.
\end{align*}
\end{example}

\begin{example}
Let $a^{ij} \in \emfb^{(ij)}[\lambda]$ and $f \in \emfb[\mu]$. Then $\widetilde{\op_{0}(a, f)} = \hnabla_{P}\hnabla_{Q}(\tilde{f}\tilde{a})^{PQ} = \tilde{a}^{PQ}\hnabla_{P}\hnabla_{Q}\tilde{f}$. Straightforward computations yield
\begin{align*}
&\hnabla_{i}\hnabla_{j}\tilde{f} = \widetilde{\nabla_{i}\nabla_{j}f - \mu P_{ij}f} = \widetilde{L_{ij}f},\\
\end{align*}
and from this and \eqref{iinfmix} there follows
\begin{align}\label{2ndorder}
\op_{0}(a, f) = a^{pq}\nabla_{p}\nabla_{q}f + \tfrac{2(1 - \mu)}{(\lambda + n+3)}\nabla_{p}a^{pq}\nabla_{q}f + \tfrac{\mu_{(2)}}{(\lambda + n + 3)_{(2)}}f\nabla_{p}\nabla_{q}a^{pq} - \tfrac{\mu(\lambda + \mu + n + 1)}{(\lambda + n + 2)}a^{pq}P_{pq}f,
\end{align}
recovering Theorem 5.3 of \cite{Bouarroudj}. When $\mu = 1$, the operator $L_{(ij)}f$ is projectively invariant, and $\op_{0}(a, f) = a^{ij}L_{(ij)}f$. 
\end{example}

\begin{example}
Let $a^{ijk} \in \emfb^{(ijk)}[\lambda]$ and $f \in \emfb[\mu]$. Then $\widetilde{\op_{0}(a, f)} = \tilde{a}^{PQR}\hnabla_{P}\hnabla_{Q}\hnabla_{R}\tilde{f}$. Straightforward computations yield
\begin{align*}
&\hnabla_{i}\hnabla_{j}\hnabla_{k}\tilde{f} = \widetilde{\nabla_{i}\nabla_{j}\nabla_{k}f} - \mu \widetilde{f\nabla_{i}P_{jk}} + 2(1 - \mu)\widetilde{P_{i(j}\nabla_{k)}f} - \mu\widetilde{P_{jk}\nabla_{i}f} = \widetilde{L_{ijk}f},
\end{align*}
and from this and \eqref{iinfmix} there follows straightforwardly
\begin{align*}
&\op_{0}(a, f) = a^{ijk}\nabla_{i}\nabla_{j}\nabla_{k}f + \tfrac{3(2-\mu)}{(\lambda + n + 5)}(\nabla_{p}a^{ijp})\nabla_{i}\nabla_{j}f\\
&\notag + ((2 - 3\mu) - \tfrac{3(\mu-1)_{(2)}}{(\lambda + n + 4)})a^{ijk}P_{ij}\nabla_{k}f + \tfrac{3(\mu-1)_{(2)}}{(\lambda + n +5)_{(2)}}(\nabla_{p}\nabla_{q}a^{ipq})\nabla_{i}f\\
&\notag +\left(-\tfrac{\mu_{(3)}}{(\lambda + n + 5)_{(3)}}\nabla_{i}\nabla_{j}\nabla_{k}a^{ijk}
+\mu(\tfrac{(\mu-1)_{(2)}}{(\lambda + n + 4)_{(2)}} - 1)a^{ijk}\nabla_{i}P_{jk} + \tfrac{\mu(\mu - 2)}{(\lambda + n + 5)}(3\mu + \tfrac{(\mu-1)}{(\lambda + n + 4)})P_{ij}\nabla_{p}a^{ijp}\right)f
\end{align*}
recovering Theorem 2.1 of \cite{Bouarroudj-third}. When $\mu = 2$ the operator $L_{(ijk)}f$ is projectively invariant and $\op_{0}(a, f) = a^{ijk}L_{(ijk)}f$.
\end{example}

\begin{example}
Let $a^{i} \in \emfb^{i}[\mu]$ and $b^{i} \in \emfb^{i}[\lambda]$. The differential operator $\op_{0}:\emf^{i}[\mu] \times \emf^{i}[\lambda] \to \emf[\mu + \lambda]$ is computed by
\begin{align*}
&\widetilde{\op_{0}(a, b)} = \hnabla_{P}\hnabla_{Q}(\tilde{a}\odot \tilde{b})^{PQ} = (\hnabla_{P}\tilde{a}^{Q})(\hnabla_{Q}\tilde{b}^{P}) + \tfrac{1}{2}\left(\tilde{a}^{P}\hnabla_{Q}\hnabla_{P}\tilde{b}^{Q}  + \tilde{b}^{P}\hnabla_{P}\hnabla_{Q}\tilde{a}^{Q} \right),\\
&\qquad\qquad\qquad\qquad\qquad\qquad\qquad = (\hnabla_{P}\tilde{a}^{Q})(\hnabla_{Q}\tilde{b}^{P}),
\end{align*}
where $\hnabla_{P}\hnabla_{Q}\tilde{b}^{P} = 0$ follows from $\hnabla_{P}\tilde{b}^{P} = 0$ and the fact that $\hnabla$ is Ricci flat. This yields
\begin{align}\label{op011}
&\op_{0}(a, b) = (\nabla_{p}a^{q})(\nabla_{q}b^{p}) - \tfrac{\mu + 1}{(n+ 1 + \lambda)}a^{p}\nabla_{p}\nabla_{q}b^{q} - \tfrac{\lambda+1}{(n+ 1 + \mu)}b^{p}\nabla_{p}\nabla_{q}a^{q}\\\notag& + \tfrac{\mu\lambda - n - 1}{(n+1+\mu)(n+1+\lambda)}(\nabla_{p}a^{p})(\nabla_{q}b^{q}) + (\mu + \lambda + 2)a^{(p}b^{q)}P_{pq}.
\end{align}
It is interesting to note the following. The differential operator $K:\emf^{i}[\lambda] \to \emf_{i}^{j}[\lambda]$ defined by $b^{i} \to K(b)_{i}^{j} = \nabla_{i}b^{j} - \tfrac{1}{n}\delta_{i}\,^{j}(\nabla_{p}b^{p})$ is invariant if and only if $\lambda = -1$. For $\lambda = \mu = -1$ it is straightforward to check that $\op_{0}(a, b)= K(a)_{p}^{q}K(b)_{q}^{p}$ and $2\op_{1}(a, b)^{i} = b^{p}K(a)_{p}^{i} + a^{p}K(b)_{p}^{i}$.
\end{example}

\begin{proposition}\label{droptwoprop}
If $a \in \emfb^{(i_{1}\dots i_{k+1})}[\lambda]$ for a non-excluded weight $\lambda$, and $u \in \emfb[k]$, then $\op_{0}(a, u) = a^{i_{1}\dots i_{k+1}}L_{i_{1}\dots i_{k+1}}u$ for a projectively invariant differential operator, $L:\emf[k] \to \emf_{(i_{1}\dots i_{k+1})}[k]$. For each choice of projective scale, $L$ has the form
\begin{align}\label{lnok}
L_{i_{1}\dots i_{k+1}}u = \nabla_{(i_{1}}\dots \nabla_{i_{k+1})}u + \sum_{s = 0}^{k-1}C^{s}_{(i_{1}\dots i_{s + 2}}\nabla_{i_{s+3}}\dots \nabla_{i_{k+1})}u, 
\end{align}
where $C^{s}_{i_{1}\dots i_{s+2}}$ is a polynomial in $P_{ij}$ and its completely symmetrized covariant derivatives of order at most $s$. Moreover there is a bijective linear map associating to each solution, $u \in \emfb[k]$, of the equation $L_{i_{1}\dots i_{k+1}}u = 0$ a homogeneity $k$ tensor $\beta \in S^{k}(T^{\ast}L)$ satisfying $\hnabla_{(I_{1}}\beta_{I_{2}\dots I_{k+1})} = 0$ and $\eul^{I_{i_{1}}}\dots \eul^{I_{i_{s}}}\hnabla_{J_{1}}\dots \hnabla_{J_{s}}\beta_{I_{1}\dots I_{k}} = 0$ for every $\{i_{1}, \dots, i_{s}\} \subset \{1, \dots l\}$.
\end{proposition}
\begin{proof}
For $u \in \emfb[\lambda]$, and a choice of projective scale define $\widetilde{L_{i_{1}\dots i_{k+1}}u} = \hnabla_{(i_{1}}\dots \hnabla_{i_{k+1})}\tilde{u}$. In general this depends on the choice of projective scale. By \eqref{iinfmix}, $\eul^{I_{i_{1}}}\dots \eul^{I_{i_{s}}}\hnabla_{I_{1}}\dots \hnabla_{I_{k+1}}\tilde{u} = 0$ for every non-empty subset $\{i_{1}, \dots, i_{s}\} \subset \{1, \dots, k+1\}$ if and only if $\lambda = k$, so that if $u \in \emfb[k]$, then $L_{i_{1}\dots i_{k+1}}u$ is a projectively invariant differential operator. In this case, by definition of $\op_{\beta}$, $\widetilde{\op_{0}(a, u)} = \tilde{a}^{I_{1}\dots I_{k+1}}\hnabla_{(I_{1}}\dots \hnabla_{I_{k+1})}\tilde{u} = \tilde{a}^{i_{1}\dots i_{k+1}}\hnabla_{(i_{1}}\dots \hnabla_{i_{k+1})}\tilde{u}$, so that $\op_{0}(a, u) = a^{i_{1}\dots i_{k+1}}L_{i_{1}\dots i_{k+1}}u$. It is straighforward to check that for $u \in \emfb[\lambda]$
\begin{align*}
L_{i_{1}\dots i_{k+1}}u = \nabla_{(i_{1}}L_{i_{2}\dots i_{k+1})}u + k(k - 1 - \lambda)P_{(i_{1}i_{2}}L_{i_{3}\dots i_{k+1})}u,
\end{align*}
and with this and $L_{ij}u = \nabla_{i}\nabla_{j}u -\lambda P_{ij}u$, an obvious inductive argument shows \eqref{lnok}.

If $u \in \emf[k]$ solves $L_{i_{1}\dots i_{k+1}}u = 0$, let $\beta_{I_{1}\dots I_{k}} = \hnabla_{(I_{1}}\dots \hnabla_{I_{k})}\tilde{u}$. By assumption and the preceeding discussion, $\hnabla_{(I_{1}}\beta_{I_{2}\dots I_{k+1})} = \hnabla_{(I_{1}}\dots \hnabla_{I_{k+1})}\tilde{u} = 0$. Similarly, \eqref{iinfmix} implies 
\begin{align*}
\eul^{I_{i_{1}}}\dots \eul^{I_{i_{s}}}\hnabla_{J_{1}}\dots \hnabla_{J_{s}}\beta_{I_{1}\dots I_{k}} = \eul^{I_{i_{1}}}\dots \eul^{I_{i_{s}}}\hnabla_{J_{1}}\dots \hnabla_{J_{s}}\hnabla_{(I_{1}}\dots \hnabla_{I_{k})} = 0. 
\end{align*}
Conversely, given $\beta$ of homogeneity $k$ satisfying the given conditions, define $u \in \emfb[k]$ by $k!\tilde{u} = \eul^{I_{1}}\dots \eul^{I_{k}}\beta_{I_{1}\dots I_{k}}$. Differentiating $\tilde{u}$ repeatedly yields $\hnabla_{(I_{1}}\dots \hnabla_{I_{k})}\tilde{u} = \beta_{I_{1}\dots I_{k}}$ and $\hnabla_{I_{1}}\dots \hnabla_{I_{k+1}}\tilde{u} = 0$, so that $L_{i_{1}\dots i_{k+1}}u = 0$.
\end{proof}

At least the existence of $L$ must have been known to the authors of \cite{Bailey-Eastwood-Gover}.

\begin{corollary}\label{tensorlemma}
If $u \in \emf[k]$ and $v \in \emf[l]$ are solutions of the projectively invariant equations $L_{i_{1}\dots i_{k+1}}u = 0$ and $L_{i_{1}\dots i_{l+1}}v = 0$, then $w = uv \in \emf[k+l]$ is a solution of $L_{i_{1}\dots i_{k+l+1}}w = 0$.
\end{corollary}

\begin{proof}
By assumption $\hnabla_{(i_{1}}\dots \hnabla_{i_{k+1})}\tilde{u} = 0$ and $\hnabla_{(i_{1}}\dots \hnabla_{i_{l+1})}\tilde{v} = 0$. As $\widetilde{uv} = \tilde{u}\tilde{v}$, each term of $\hnabla_{(i_{1}}\dots \hnabla_{i_{k+l+1})}\widetilde{uv}$ involves the complete symmetrization of either at least $k+1$ derivatives of $\tilde{u}$ or at least $l+1$ derivatives of $\tilde{v}$, so must vanish.
\end{proof}

Corollary \ref{tensorlemma} reflects that in the flat model case the space of solutions of $L_{i_{1}\dots i_{k+1}}u = 0$ is canonically identified with the $k$th symmetric power of the standard representation, $\op$, of $\mathfrak{sl}(n+1, \rea)$. 

\begin{proposition}\label{ricflatprop}
Suppose $a \in \emfb^{(i_{1}\dots i_{k})}[\lambda]$ ($\lambda$ is not an excluded weight) and $f \in \emfb[\mu]$. If $\nabla \in \en$ is a Ricci flat representative then on $M$ there holds
\begin{align}\label{ricflatlhf}
\op_{\be}(a, f) = \sum_{m = 0}^{k}\frac{\binom{k-\be-1-\mu}{m}\binom{k-\be}{m}}{\binom{n+2k -1 + \lambda}{m}}(\nabla_{j_{1}}\dots \nabla_{j_{m}}a^{i_{1}\dots i_{k-\be - m}j_{1}\dots j_{m}})(\nabla_{i_{1}}\dots \nabla_{i_{k-\be - m}}f).
\end{align}
 For a general $\en$, if $\nabla \in \en$ is the representative associated to a projective normal scale at $x \in M$ then at the point $x$ there holds \eqref{ricflatlhf}.
\end{proposition}
\begin{proof}
In the case that $\en$ admits a Ricci flat representative, using \eqref{ricflatlift} and \eqref{iinfmix} it is not hard to show that \eqref{ricflatlhf} holds everywhere. For the most important case, $\beta = 0$, this was shown by M. Bordemann, \cite{Bordemann}. From the definition of $\op_{\be}$ and $\hnabla$ and the reasoning of the proof of Proposition \ref{droptwoprop} it is apparent that in general the explicit expression for $\op_{\be}(a, f)$ can differ from that holding in the Ricci flat case only by terms of order at most $k-\be-2$ in $f$. Moreover, for $0 \leq s \leq k-\be-2$, the coefficient of a term containing $\nabla_{i_{1}}\dots \nabla_{i_{s}}f$ is a polynomial in $a$ and its covariant derivatives of order at most $k-\be-s$ and in $P_{ij}$ and its covariant derivatives of order at most $k-\be-2-s$. Every appearance of $\nabla_{i_{1}\dots i_{r}}P_{jk}$ is contracted with expressions involving $a$ and its covariant derivatives, all of which are completely symmetric, and so every appearance of $\nabla_{i_{1}\dots i_{r}}P_{jk}$ in the expression for $\op_{\be}(a, f)$ can be replaced by $\nabla_{(i_{1}\dots i_{r}}P_{jk)}$. In a projective normal scale at $x$, all such expressions vanish at $x$, and this proves the claim.
\end{proof}
\begin{remark}\label{flatcaseremark}
In the flat case, when $\lambda = 0$ and $\beta = 0$, \eqref{ricflatlhf} recovers the projectively invariant quantization map defined in equations 4.13 and 4.14 of \cite{Lecomte-Ovsienko}. See Section \ref{quantizationsection} below for further discussion.
\end{remark}

The following proposition shows that the operators $\op_{\beta}$ arise naturally as a measure of how the invariant lift fails to preserve the symmetric product.

\begin{proposition}\label{liftsymprop}
  For $1 \leq q \leq p$ let $a_{q} \in \emfb^{(i_{1}\dots i_{k_{q}})}[\lambda_{q}]$ for non-excluded weights $\lambda_{q}$ such that for every $J \in \{1, \dots, p\}$, the sums $\sum_{q \in J}\lambda_{j}$ are also non-excluded. Write $K = \sum_{q = 1}^{p}k_{q}$ and $\Lambda = \sum_{q = 1}^{p}\lambda_{q}$. Then
\begin{align}\label{liftsym}
 &\widetilde{a_{1}\odot \dots \odot a_{p}}^{I_{1}\dots I_{K}} =\\ 
& \notag (\tilde{a}_{1} \odot \dots \odot \tilde{a}_{p})^{I_{1}\dots I_{K}} - \sum_{s = 1}^{K} \frac{\binom{K}{s}}{(\Lambda + n + 2K - s)_{(s)}}\eul^{(I_{1}}\dots \eul^{I_{s}}\widetilde{\op_{K-s}(a_{1}, \dots, a_{p})}\,^{I_{s+1}\dots I_{K})}.
\end{align}
\end{proposition}
\begin{proof}
Suppose that there exist $u_{r} \in \emfb^{(i_{1}\dots i_{k-r})}[\Lambda]$ such that
\begin{align*}
\widetilde{a_{1}\odot \dots \odot a_{p}}^{I_{1}\dots I_{K}} - \widetilde{a}_{1}\odot \dots \odot \tilde{a}_{p}^{I_{1}\dots I_{K}} = -\sum_{s = 1}^{K}\eul^{(I_{1}}\dots \eul^{I_{s}}\tilde{u}_{s}^{I_{s+1}\dots I_{K})}.
\end{align*}
If $v^{I_{1}\dots I_{k}} = \eul^{(I_{1}}\dots \eul^{I_{s}}\tilde{b}^{I_{s+1}\dots I_{k})}$ for $b \in \emfb^{(i_{1}\dots i_{k-s})}[\lambda]$, then it can be checked that 
\begin{align*}
\hnabla_{P}v^{I_{1}\dots I_{k-1}P} = \tfrac{s(\lambda + n + 2k -s)}{k}\eul^{(I_{1}}\dots \eul^{I_{s-1}}\tilde{b}^{I_{s}\dots I_{k-1})},
\end{align*}
and applying this observation with $b = u_{s}$ and using the defintion of $\widetilde{a_{1}\odot \dots \odot a_{p}}$ yields
\begin{align}\label{liftsym2}
\hnabla_{I_{1}}\dots \hnabla_{I_{r}}(\tilde{a}_{1}\odot \dots \odot \tilde{a}_{p})^{I_{1}\dots I_{K}} = \sum_{s = r}^{K}\frac{s_{(r)}(\Lambda + 2K + n -s)_{(r)}}{K_{(r)}}\eul^{(I_{1}}\dots \eul^{I_{s-r}}\tilde{u}_{s}^{I_{s-r+1}\dots I_{K-r})}.
\end{align}
The $s = r$ term of this sum and the definition of $\op_{K-r}$ yield 
\begin{align*}
u_{r}^{i_{1}\dots i_{K-r}} = \frac{\binom{K}{r}}{(\Lambda + 2K + n - r)_{(r)}}\op_{K-r}(a_{1}, \dots, a_{p})^{i_{1}\dots i_{K-r}},
\end{align*}
and this proves the claim.
\end{proof}

\subsection{One-Dimensional Case}\label{onedimsection}
The one-dimensional case of Theorem \ref{maintheorem} is more interesting than it may appear at first. With the canonical identification $\emf^{(i_{1}\dots i_{k})}[\lambda] \simeq \emf[\lambda + 2k]$, Proposition \ref{invliftprop} implies that a section of $\emf[\mu]$ admits an invariant lift to a symmetric $k$-tensor as long as $\mu$ is not in the set of excluded weights $\{0,1, \dots, k-1\}$. For a non-negative integer $k$, $\sigma \in \rea  - \{0,1, \dots, k-1\}$, and $u \in \emfb[\sigma]$, define $\tk(u)$ to be the invariant lift of $u$ viewed as a section of $\emf^{(i_{1}\dots i_{k})}[\sigma - 2k]$, so that
\begin{align}
\label{onedlift}&(-1)^{k}\binom{\sigma}{k}\tk(u) = (-1)^{k}\binom{\sigma}{k}\sum_{m = 0}^{k}\frac{(-1)^{m}\binom{k}{m}}{(\sigma)_{(m)}}\widetilde{D^{m}(u)}\hat{E}^{k-m} \odot \eul^{m}\\
\notag & = \sum_{m = 0}^{k}\binom{k-1-\sigma}{k-m}\frac{1}{m!}\widetilde{D^{m}(u)}\hat{E}^{k-m} \odot \eul^{m},
\end{align}
Here $D$ is the divergence operator determined by a representative of a flat projective structure in one-dimension, and $E$ is a unimodular frame corresponding to a volume form made parallel by $D$. If $1 \leq m \leq p$, $\sigma_{m} \notin \{0, \dots, k_{m}-1\}$, $K = \sum_{m = 1}^{p}k_{m}$, and $\Sig = \sum_{m = 1}^{p}\sigma_{m}$, the operator $\op_{\be}:\emfb[\lambda_{1} + 2k_{1}]\times \dots \emfb[\lambda_{p} + 2k_{p}] \to \emfb[\Lambda + 2\be]$ determines a multilinear pairing,
\begin{align*}
&\rc_{\be}^{k_{1},\dots, k_{p}}:\emfb[\sigma_{1}] \times \dots \times \emfb[\sigma_{p}] \to \emfb[\Sig - 2(K - \be)],& &\text{defined by}\\
&\rc_{\be}^{k_{1}, \dots, k_{p}}(u_{1}, \dots, u_{p}) = \op_{\be}(\tk_{1}(u_{1}), \dots, \tk_{p}(u_{p})).
\end{align*}
Straightforward computations using \eqref{onedlift} and $(\tr \hnabla)^{l} (\hat{E}^{k-p}\odot \eul^{p}) = \tfrac{(p)_{(l)}(2k+1-p)_{(l)}}{k_{(l)}}\hat{E}^{k-p}\odot \eul^{p-l}$ yield
\begin{align*}
&\notag&&(-1)^{K}(\prod_{i = 1}^{p}\binom{\sigma_{i}}{k_{i}})\rc^{k_{1}, \dots, k_{p}}_{0}(u_{1}, \dots, u_{p})  \\
&\notag =&& (\tr \hnabla)^{K}\left \{ \sum_{m = 0}^{k}(\prod_{i = 1}^{p}\binom{k_{i} - 1 - \sigma_{i}}{k_{i} - m_{i}})\prod_{i = 1}^{p}\tfrac{D^{m_{i}}(u_{i})}{m_{i}!}) \hat{E}^{K - M}\odot \eul^{M}  \right \}  \\
&\notag =&& \sum_{m = 0}^{k}\sum_{s = 0}^{M}(\prod_{i = 1}^{p}\tfrac{\binom{k_{i} - 1 - \sigma_{i}}{k_{i} - m_{i}}}{m_{i}!})\tbinom{K}{s} (\tr)^{K}\left \{ \hnabla^{K-s}(\prod_{i=1}^{p}D^{m_{i}}(u_{i}))\tensor \hnabla^{s}(\hat{E}^{K-M}\odot \eul^{M}) \right \}  \\
& \notag =&& \sum_{m = 0}^{k}\sum_{s = 0}^{M}(\prod_{i = 1}^{p}\tfrac{\binom{k_{i} - 1 - \sigma_{i}}{k_{i} - m_{i}}}{m_{i}!})\tbinom{M}{s}\tbinom{2K+1-M}{s}s! (\tr)^{K-s}\left \{ \hnabla^{K-s}(\prod_{i=1}^{p}D^{m_{i}}(u_{i}))\tensor \hnabla^{s}(\hat{E}^{K-M}\odot \eul^{M-s}) \right \}  \\
& \notag =&& \sum_{m= 0}^{k}\sum_{s = 0}^{M}(\prod_{i = 1}^{p}\tbinom{k_{i} - 1 - \sigma_{i}}{k_{i} - m_{i}})\tbinom{M}{m_{1}; \dots ; m_{p}}\tbinom{2K+1-M}{s}\tbinom{\Sigma - 3K+M}{M-s} D^{K-M}\left (\prod_{i=1}^{p}D^{m_{i}}(u_{i})) \right ),  \\
\end{align*}
where $M = \sum_{i = 1}^{p}m_{i}$, $\binom{M}{m_{1};\dots;m_{p}}$ is the standard multinomial coefficient, and $\sum_{m = 0}^{k}$ is written sometimes for $ \sum_{m_{1} = 0}^{k_{1}}\dots \sum_{m_{p} = 0}^{k_{p}}$. Simplifying yields
\begin{align}
&\label{rckp}&&(-1)^{K}(\prod_{i = 1}^{p}\binom{\sigma_{i}}{k_{i}})\rc^{k_{1}, \dots, k_{p}}_{0}(u_{1}, \dots, u_{p})  \\
&\notag =&&\sum_{m_{1} = 0}^{k_{1}}\dots \sum_{m_{p} = 0}^{k_{p}}(\prod_{i = 1}^{p}\binom{k_{i} - 1 - \sigma_{i}}{k_{i} - m_{i}})\binom{M}{m_{1};\dots ;m_{p}}\binom{\Sig - K + 1}{M}D^{K-M}(\prod_{i = 1}^{p}D^{m_{i}}(u_{i})),
\end{align}
The given form of \eqref{rckp} is not satisfactory; it admits many simplifications and it is not clear which is the most useful. The identities $\binom{\sigma}{m}\binom{k- 1 - \sigma}{k-m} = (-1)^{k-m}\binom{\sigma}{m}\binom{\sigma -m}{k - m} = (-1)^{k-m} \binom{\sigma}{k}\binom{k}{m}$ suggest some manipulations. 

Up to a normalizing constant, the \textbf{Rankin-Cohen brackets}, $\rc_{k}:\emfb[\sigma_{1}]\times \emfb[\sigma_{2}] \to \emfb[\sigma_{1} + \sigma_{2} -2k]$, are defined (see Theorem 7.1 of \cite{Cohen-lift} or eq. 3.1 of \cite{Cohen-Manin-Zagier}) by
\begin{align*}
\rc_{k}(u_{1}, u_{2}) = \sum_{m = 0}^{k}(-1)^{m}\binom{k-1-\sigma_{2}}{m}\binom{k-1-\sigma_{1}}{k-m}D^{m}(u_{1})D^{k-m}(u_{2}).
\end{align*}
(The correspondence with usual notations is given by regarding a modular form of weight $k$ as a section of $\emfb[-k]$). The Rankin-Cohen brackets are graded skew-symmetric in the sense that $\rc_{k}(u_{2}, u_{1}) = (-1)^{k}\rc_{k}(u_{1}, u_{2})$. The notation for $\rc_{k}$ ought to include the weights $\sigma_{1}$ and $\sigma_{2}$, e.g. $\rc_{k, \sigma_{1}, \sigma_{2}}$, and that with such an extended notation skew-symmetry reads $\rc_{k, \sigma_{1}, \sigma_{2}}(u_{1}, u_{2}) = (-1)^{k}\rc_{k, \sigma_{2}, \sigma_{1}}(u_{2}, u_{1})$. 

Expanding \eqref{rckp} using the Leibniz rule and the identity $\sum_{s = 0}^{p}\binom{a}{s}\binom{b}{p - s} = \binom{a + b}{p}$ shows that 
\begin{align*}
&(-1)^{k}\binom{\sigma_{1}}{k}\rc^{k, 0, \dots, 0}_{0}(u_{1}, \dots, u_{p}) 
= \sum_{m = 0}^{k}\sum_{q = m}^{k}\binom{k - 1- \si_{1}}{k- m}\binom{\Sigma - k +1}{m}\binom{k-m}{q - m}D^{q}(u_{1})D^{k-q}(\prod_{i = 2}^{p}u_{i})\\
& =  \sum_{q = 0}^{k}\sum_{m = 0}^{q}\binom{k - 1- \si_{1}}{k- m}\binom{\Sigma - k +1}{m}\binom{k-m}{q - m}D^{q}(u_{1})D^{k-q}(\prod_{i = 2}^{p}u_{i})\\
& =  \sum_{q = 0}^{k}\sum_{m = 0}^{q}\binom{k - 1- \si_{1}}{k- q}\binom{\Sigma - k +1}{m}\binom{q-\si_{1}}{q - m}D^{q}(u_{1})D^{k-q}(\prod_{i = 2}^{p}u_{i})\\
& =  \sum_{q = 0}^{k}(-1)^{q}\binom{k - 1- \si_{1}}{k- q}\binom{k - 1- (\Sigma - \si_{1})}{q}D^{q}(u_{1})D^{k-q}(\prod_{i = 2}^{p}u_{i}) = \rc_{k}(u_{1}, \prod_{i = 2}^{p}u_{i}).
\end{align*}
In particular, $\rc_{k}(u, v) = (-1)^{k}\binom{\sigma_{1}}{k}\rc_{0}^{k, 0}(u, v)$. This is the sense in which the operators $\op_{\be}$ generalize the Rankin-Cohen brackets.

Similarly, direct computation using \eqref{onedlift} (or \eqref{ricflatlhf}) shows that if $\lambda +2k \notin \{0, \dots, k-1\}$ then $\op_{\be}:\emfb[\lambda + 2k]\times \emfb[\mu] \to \emfb[\lambda + \mu + 2\be]$ is given by
\begin{align*}
\op_{\be}(u, v) = \sum_{m = 0}^{k-\be}\frac{\binom{k-\be - 1-\mu}{m}\binom{k-\be}{m}}{\binom{2k + \lambda}{m}}D^{m}(u)D^{k-\beta-m}(v) = \sum_{m = 0}^{k-\be}\frac{\binom{k-\be - 1-\mu}{m}\binom{k}{m}\binom{k-m}{\be}}{\binom{2k + \lambda}{m}\binom{k}{\be}}D^{m}(u)D^{k-\beta-m}(v),
\end{align*} 
which recovers a special case of equation 5.3 of Proposition 7 of \cite{Cohen-Manin-Zagier}. The associated pairing, $\rc_{\be}^{k, 0}:\emfb[\sigma_{1}]\times \emfb[\sigma_{2}] \to \emfb[\sigma_{1} + \sigma_{2} -2(k-\be)]$, determined by $\op_{\be}$ (and defined for $\sigma_{i} \notin \{0, \dots, k-1\}$) is 
\begin{align*}
&(-1)^{k}\binom{\si_{1}}{k}\binom{k}{\be}\rc_{\be}^{k, 0}(u_{1}, u_{2}) \\
&= \sum_{m = 0}^{k-\be}(-1)^{m}\binom{k-\be-1-\sigma_{2}}{m}\binom{k-1-\si_{1}}{k - m}\binom{k-m}{\be}D^{m}(u_{1})D^{k-\be-m}(u_{2})\\
& = \binom{k - 1- \si_{1}}{\be}\sum_{m = 0}^{k-\be}(-1)^{m}\binom{k-\be - 1-\sigma_{2}}{m}\binom{k-\be - 1-\si_{1}}{k -\be - m}D^{m}(u_{1})D^{k-\be-m}(u_{2})\\
& = (-1)^{k-\beta}\binom{\si_{1}}{k-\be}\binom{k - 1- \si_{1}}{\be}\rc^{k-\be, 0}_{0}(u_{1}, u_{2}) = (-1)^{k}\binom{\si_{1}}{k}\binom{k}{\be}\rc^{k-\be, 0}_{0},
\end{align*}
which shows that $\rc^{k, 0}_{\be} = \rc^{k-\be, 0}_{0}$. %(Of course this follows also from the uniqueness of the Rankin-Cohen brackets).

\subsection{Another Definition of $\op_{\beta}$}\label{anothersection}
For any smooth manifold, $N$, the action of $\diff(N)$ on $T^{\ast}N$ is Hamiltonian with respect to the tautological Poisson structure, $\pois$. A moment map, $\Phi:\vect(N) \to C^{\infty}(T^{\ast}N)$, for this action is defined as follows. Each vector field, $X \in \vect(N)$ has a tautological lift to a vector field $\bar{X} \in \vect(T^{\ast}N)$ which is the infinitesimal generator of the action on $T^{\ast}N$ of the differential of the flow of $X$. The moment map, $\Phi$, is defined by $\Phi(X) = \alpha(\bar{X})$, where $\alpha$ is the tautological one-form on $T^{\ast}N$. $\Phi$ extends to give the tautological algebra isomorphism $\Phi:\sy(N) \to \pol(T^{\ast}N)$, and this extension is $\diff(N)$-equivariant in the sense that $\Phi(\lie_{X}a) = \lie_{\bar{X}}\Phi(a)$. Let $\teul$ denote the Euler vector field generating the dilations in the fibers of $T^{\ast}N$. If $\nabla$ is an affine connection on $N$, and if $a \in \Gamma(S^{k}(TN))$, denote by $\sdiv(a)$ the section of $S^{k}(TN)$ defined by $\nabla_{p}a^{i_{1}\dots i_{k-1}p}$.

\begin{lemma}\label{divlemma}
If $N$ is a smooth manifold equipped with an affine connection, $\nabla$, there exists a unique linear second-order differential operator $D:C^{\infty}(T^{\ast}N) \to C^{\infty}(T^{\ast}N)$ such that for each $k \in \nat$, 
\begin{align}\label{kdiv}
&D(\Phi(a)) = k\Phi(\sdiv(a)),& \text{for all $a \in \Gamma(S^{k}(TN))$}.
\end{align}
\end{lemma}
\begin{proof}
If $D^{\prime}$ is a second linear second-order differential operator satisfying the given conditions, then $E = D^{\prime} - D$ is a linear differential operator of order at most $2$ annihilating $\pol(T^{\ast}N)$. For $f \in C^{\infty}(N)$ and $a \in \Gamma(S^{k}(TN))$, 
\begin{align*}
D(\Phi(fa)) = k\Phi(\sdiv(fa)) = k\Phi(f)\Phi(\sdiv(a)) + k\Phi(df(a))= \Phi(f)D(\Phi(a)) + k\Phi(df(a)),
\end{align*} 
and so $E(\Phi(fa)) = \Phi(f)E(\Phi(a))$, so that the differential operator $\Phi^{-1} \circ E\circ \Phi$ is a $C^{\infty}(N)$-module map from $\pol(T^{\ast}N) \to \pol(T^{\ast}N)$. Because this map vanishes on each graded piece $\Gamma(S^{k}(TN))$ it must be identically $0$. A linear second order differential operator annihilating $\pol(T^{\ast}N)$ must vanish on $C^{\infty}(T^{\ast}N)$, and so $E = 0$. This shows uniqueness.

If $x^{i}$ are local coordinates on $N$, let $z_{i}$ denote the natural coordinates on the fiber $T_{x}^{\ast}N$ (so that $\alpha = z_{i}dx^{i}$ and $\teul = z_{i}\tfrac{\partial}{\partial z_{i}}$), and let $\Gamma_{ij}\,^{k}$ denote the Christoffel symbols of $\nabla$ with respect to the coordinate frame $\tfrac{\pr}{\pr x^{i}}$. Keeping in mind that the components, $\Gamma_{ij}\,^{k}$, do not transform tensorially, it is straightforward to check that $D = \tfrac{\partial^{2}}{\partial z_{i}\partial x^{i}} + z_{k}\Gamma_{ij}\,^{k}\tfrac{\pr^{2}}{\pr z_{i}\pr z_{j}} + \Gamma_{ip}\,^{p}\tfrac{\partial}{\partial z_{i}}$ does not depend on the choice of coordinates and is a differential operator satisfying the given conditions. 
\end{proof}
Note that the operator $D$ of Lemma \ref{divlemma} has homogeneity $-1$ in the sense that $[\lie_{\teul}, D] = -D$. 

Applying Lemma \ref{divlemma} to the ambient connection, $\hnabla$, yields a second order differential operator $\hD:C^{\infty}(T^{\ast}\form) \to C^{\infty}(T^{\ast}\form)$, and an algebra isomorphism $\hPhi:\sy(\form) \to \pol(T^{\ast}\form)$ (ordinarily in this paper the notation for the isomorphism $\hPhi$ has been suppressed). The operators $\op_{\beta}$ may be defined in terms of $\hD$ by
\begin{align*}
K_{(K-\beta)}\widetilde{\op_{\beta}(a_{1}, \dots, a_{p})} = \hPhi^{-1}(\hD^{K - \beta}(\hPhi(\tilde{a}_{1}\odot \dots \odot \tilde{a}_{p}))) \mod \eul.
\end{align*}

The entire function $\sigma(z) = \tfrac{1}{\Gamma(1 - z)} = \sum_{m \geq 0}c_{m}z^{m}$ has zeros only at the postive integers. To the differential operator, $z \tfrac{d}{dz}$, on the line, associate the formal pseudo-differential operator $\sigma(z\tfrac{d}{dz}) = \sum_{m \geq 0}c_{m}(z\tfrac{d}{dz})^{m}$ and observe that $\sigma(z \tfrac{d}{dz})(z^{k}) = \sigma(k)z^{k} = 0$ for $k$ a positive integer and that $\sigma(z \tfrac{d}{dz})(1) = c_{0} = \sigma(0) = 1$. Consequently, it makes sense to apply $\sigma(z \tfrac{d}{dz})$ to any function analytic and entire in $z$, and the result is the zeroth term in the series expansion of the function. Define $\pdop = \sum_{m \geq 0}\tfrac{1}{m!}\sigma(z \tfrac{d}{dz}) \circ (\tfrac{d}{dz})^{m}$. If $f$ is any entire function, $\pdop(f)$ is defined and equals $f(1)$. In particular, $\pdop(z^{k}) = 1$ for every non-negative integer $k$. 

Now let $\hteul$ denote the Euler vector field generating the dilations in the fibers of $T^{\ast}\form$. Define a formal pseudo-differential operator $\pdop = \sum_{m\geq 0}\tfrac{1}{m!}\sigma(\lie_{\hteul})\circ \hD^{m}$. Reasoning as in the one-dimensional case shows that, for $A \in S^{k}(T\form)$, 
\begin{align*}
&\pdop(\hPhi(A)) %= \sum_{m\geq 0}\tfrac{1}{m!}\sigma(\lie_{\hteul})\circ \hD^{m}(\hPhi(A)) 
= \sum_{m= 0}^{k}\tfrac{1}{m!}\sigma(\lie_{\hteul})\circ \hD^{m}(\hPhi(A)) = \sum_{m =  0}^{k}\tfrac{\sigma(k-m)}{m!}\hD^{m}(\hPhi(A)) = \tfrac{1}{k!}\hD^{k}(\hPhi(A)).
\end{align*}
It follows that for $a_{i} \in \alg_{k_{i}}$
\begin{align*}
\pdop(\hPhi(\tilde{a}_{1}\odot \dots \odot \tilde{a}_{p})) = \tfrac{1}{K!}\hD^{k}(\hPhi(\tilde{a}_{1}\odot \dots \odot \tilde{a}_{p})) = \hPhi(\widetilde{\op_{0}(a_{1}, \dots, a_{p})}).
\end{align*}
Alternatively, one can write
\begin{align*}
\widetilde{\op_{0}(a_{1}, \dots, a_{p})} = \sum_{m \geq 0}\hPhi^{-1}\circ\sigma(\lie_{\hteul})\circ \hPhi \circ (\sdiv)^{m}(\tilde{a}_{1}\odot \dots \odot \tilde{a}_{p}).
\end{align*}

\section{Star Products}\label{spsection}
\subsection{Projectively Invariant Deformation Quantization}\label{quantizationsection}
Let $\dop^{k}_{\mu, \lambda} \subset \dop_{\mu, \lambda}$ be as in the introduction. For $\lambda \notin \{-n, -n -1, -n -2, \dots \}$ and $\mu \notin \{-1, -2, \dots\}$ define the \textbf{projectively invariant quantization map} $\op:\alg_{\lambda} \to \dop_{\mu, \lambda}$ by setting $\op(a)(u) = \op_{0}(a, u)$ for $a \in \alg_{\lambda, k}$ and $u \in \emfb[\mu]$, and extending linearly to $\alg_{\lambda}$. To any $A \in \alg_{\lambda}$ there corresponds an integer $k \geq 0$, and a sequence of $a_{\ga} \in \alg_{\lambda, k}$ for $1 \leq \ga \leq k$, such that $A = \sum_{\ga = 0}^{k} a_{\ga}$ and hence an $\op(A) \in \dop_{\mu, \lambda}^{k}$ defined by $\op(A) = \sum_{\ga = 0}^{k}\op(a_{\ga})$. Because the principal symbol of $\op(A)$ is the top order piece of $A$, $\op:\alg_{\lambda} \to \dop_{\mu, \lambda}$ is an injective map. The principal symbol of $\dop \in \dop_{\mu, \lambda}^{k}$ is some $a_{k} \in \alg_{\lambda, k}$, and so $\dop - \op(a_{k}) \in \dop_{\mu, \lambda}^{k-1}$ is a differential operator of order at most $k-1$.  The principal symbol of $\dop - \op(a_{k}) \in \dop_{\mu, \lambda}^{k-1}$ is some $a_{k-1} \in \alg_{\lambda, k-1}$, and so $\dop - \op(a_{k}) - \op(a_{k-2}) \in \dop^{k-2}_{\lambda, \mu}$. Iterating this procedure yields, for $0 \leq \gamma \leq k$, a sequence of $a_{\ga} \in \alg^{\lambda, k-\ga}$ (some of which may be $0$), so that $\dop = \sum_{\gamma = 0}^{k}\op(a_{\gamma}) = \op(A)$, where $A = \sum_{\ga = 0}^{k}a_{\ga}\in \alg_{\lambda}$. This shows that $\op$ is onto and admits a linear inverse, $\op^{-1}:\dop_{\mu, \lambda}^{k} \to \alg_{\lambda}$, which should be called a \textbf{projectively invariant symbol map}, as in the projectively flat case.

 As explained in Remark \ref{flatcaseremark} equation \eqref{ricflatlhf} shows that $\op$ is the projectively invariant quantization map defined in the flat case by Lecomte-Ovsienko, \cite{Lecomte-Ovsienko}. Moreover, Proposition \ref{ricflatprop} shows that if $\nabla \in \en$ is the representative associated to a choice of projective normal scale at $x \in M$, then $\op(a)$ satisfies equation \eqref{ricflatlhf} (with $\beta = 0$ and $\lambda = 0$) at the point $x$.

\subsection{Self-Adjointness of Quantization Map and Half-Densities}
On the space, $\emf_{c}[\mu]$, of compactly supported sections of $\emf[\mu]$, there is canonical pairing, $(\,,\,):\emf_{c}[\mu] \times \emf_{c}[-n-1-\mu] \to \rea$ defined by integration, $(u, v) = \int_{M}uv$. The specialization to the flat case of the following proposition was shown in \cite{Duval-Lecomte-Ovsienko} by invoking uniqueness of the flat projectively invariant quantization map. 

\begin{theorem}\label{selfadjoint}
If $a^{i_{1}\dots i_{k}} \in \emfb^{(i_{1}\dots i_{k})}$ then $\op(a)$ is self-adjoint in the sense that $(\op(a)(u), v) = (-1)^{k}(u, \op(a)(v))$ for $u \in \emfb_{c}[\mu]$ and $v \in \emfb_{c}[-n-1-\mu]$.
\end{theorem}

\begin{proof}
First let $g \in \emfb[-n-1]$, so that $\int_{M}g$ is defined. Let $s:U \subset M \to \form$ be any choice of projective scale on the open domain $U$. It is claimed that $\int_{U}g = \int_{s(U)}\tilde{g}\alpha$, where $\alpha$ is the tautological $n$-form on $\form$. Because $\alpha$ has homogeneity $n+1$, the $n$-form $\tilde{g}\alpha$ has homogeneity $0$, and so replacing $s$ in the integral by $\tilde{s} = fs$ does not change the value of the integral. The restriction to $U$ of $g$ equals $hs^{n+1}$ for some $h \in C^{\infty}(U)$, and by definition there hold $s^{\ast}(\tilde{g}) = h$ and $s^{\ast}(\alpha) = s^{n+1}$, so that $\int_{U}g = \int_{U}hs^{n+1} = \int_{s(U)}\tilde{g}\alpha$.

If $A\in \vect(\form)$ has homogeneity $-n-1$, then $\sdiv_{\Psi}(A)$ is a function on $\form$ of homogeneity $-n-1$, and so (abusing notation slightly), it makes sense to write $\int_{M}\sdiv_{\Psi}(A)$ to denote the integral of the corresponding $n$-form on $M$. It is claimed that $\int_{M}\sdiv_{\Psi}(A) = 0$. First it is shown that for $A \in \vect(\form)$ of homogeneity $\lambda$,
\begin{align}\label{divalpha}
(n+1)\sdiv_{\Psi}(A)\alpha = (n+1+\lambda)i(A)\Psi - \lambda d(i(A)\al).
\end{align}
By definition $\sdiv_{\Psi}(A)\Psi = \lie_{A}(\Psi) = \lie_{A}(d\al) = d(\lie_{A}\al)$. Interior multiplying $\eul$ gives 
\begin{align}
\label{divalinter}&(n+1)\sdiv_{\Psi}(A)\al = i(\eul)d(\lie_{A}\al) = \lie_{\eul}(\lie_{A}\al) - d(i(\eul)\lie_{A}\al)\\
\notag & = (n+1+\lambda)\lie_{A}\al - d(i(\eul)\lie_{A}\al).
\end{align}
Similarly,
\begin{align*}
&i(\eul)\lie_{A}\al = i(\eul)(i(A)d\al + d(i(A)\al)) = -i(A)(i(\eul)d\al) + \lie_{\eul}(i(A)\al) - d(i(\eul)(i(A)\al))\\
& = -(n+1)i(A)\al + (n+1+\lambda)i(A)\al = \lambda i(A)\al,
\end{align*}
and in \eqref{divalinter} this gives \eqref{divalpha}. In particular, in the case that $\lambda = -n -1$ there holds $\sdiv_{\Psi}(A)\al = d(i(A)\al)$. Choose a locally finite open cover, $\{U_{\be}\}$, of $M$, a subordinate partitition of unity, $\{\phi_{\be}\}$, and on each $U_{\be}$ a local section $s_{\be}:U_{\be} \to \form$. Then Stokes's Theorem yields
\begin{align*}
&\int_{M}\sdiv_{\Psi}(A) = \sum_{\be}\int_{U_{\be}}\phi_{\be}\sdiv_{\Psi}(A) = \sum_{\be}\int_{s_{\be}(U_{\be})}(\phi_{\be}\circ \rho)\sdiv_{\Psi}(A) = \sum_{\be}\int_{s_{\be}(U_{\be})}(\phi_{\be}\circ \rho) d(i(A)\al) \\
& = \sum_{\be}\int_{s_{\be}(U_{\be})}d((\phi_{\be}\circ \rho) i(A)\al) - \sum_{\be}\int_{s_{\be}(U_{\be})}d(\phi_{\be}\circ \rho) \wedge i(A)\al
 = - \sum_{\be}\int_{s_{\be}(U_{\be})}d(\phi_{\be}\circ \rho) \wedge i(A)\al \\
&= - \sum_{\be}\int_{U_{\be}}d\phi_{\be} \wedge s_{\be}^{\ast}(i(A)\al) = -\int_{M}d(\sum_{\be}\phi_{\be}) \wedge s_{\be}^{\ast}(i(A)\al) = 0.
\end{align*}
so that $\int_{M}\sdiv_{\Psi}(A) = 0$.

The proposition will be proved now by an integration by parts argument which uses the preceeding observation. For $0 \leq s \leq k - 1$, define vector fields $B_{s}^{I}$ of homogeneity $-n-1$ by 
\begin{align*}
B_{s}^{I} = \tilde{a}^{II_{1}\dots I_{k-1}}(\hnabla_{I_{1}}\dots \hnabla_{I_{s}}\tilde{v})(\hnabla_{I_{s+1}}\dots \hnabla_{I_{k-1}}\tilde{u}).
\end{align*}
Note that because $\hnabla \Psi = 0$ there holds $\sdiv_{\Psi}(B) = \hnabla_{I}B^{I}$ for any $B \in \vect(\form)$. Repeatedly using together the Leibniz rule and $\tr \hnabla \tilde{a} = 0$ yields
\begin{align}\label{intpartsid}
\widetilde{\op(a)(u)}\tilde{v} = (-1)^{k}\widetilde{\op(a)(v)}\tilde{u} + \sum_{s = 0}^{k-1}(-1)^{s}\sdiv_{\Psi}(B_{s})
\end{align}
Since the $B_{s}$ have homogeneity $-n-1$ integrating over $M$ the $-(n+1)$-densities corresponding to \eqref{intpartsid} yields $(\op(a)(u), v) = (-1)^{k}(u, \op(a)(v))$.
\end{proof}
Theorem \ref{selfadjoint} shows that the representation of $\alg$ as operators on the space of \textbf{half-densities} (sections of $\emf[-(n+1)/2]$) is formally self-adjoint.

\subsection{Projectively Invariant Star Product}
On a smooth manifold, $N$, restricting to $\pol(T^{\ast}N)$ the tautological Poisson bracket on $C^{\infty}(T^{\ast}N)$ and using the identification of $\pol(T^{\ast}N)$ with the algebra of symmetric tensors, $\sy(N)$, yields the symmetric Schouten bracket, $\{\,,\,\}:\sy(N) \times \sy(N) \to \sy(N)$. If $\nabla$ is any torsion-free affine connection on $N$ set
\begin{align}\label{schouten}
\{A, B\}^{i_{1}\dots i_{k+l-1}} = kA^{p(i_{1}\dots i_{k-1}}\nabla_{p}B^{i_{k}\dots i_{k+l-1})} - lB^{p(i_{1}\dots i_{l-1}}\nabla_{p}A^{i_{l}\dots i_{k+l-1})}.
\end{align}
Using the identity 
\begin{align}
A^{(pi_{1}\dots i_{k-1}}B^{i_{k}\dots i_{k+l-1})} = \tfrac{k}{k+l}A^{p(i_{1}\dots i_{k-1}}B^{i_{k}\dots i_{k+l-1})} + \tfrac{l}{k+l}B^{p(i_{1}\dots i_{l-1}}A^{i_{l}\dots i_{k+l-1})},
\end{align}
it is easy to check directly the Leibniz rule, $\{A\odot B, C\} = A\odot \{B, C\} + B \odot \{A, C\}$. It is straightforward to check directly that for $A, B, C \in \vect(N)$, there holds $\text{Cycle} \{\{A, B\}, C\} = 0$, and from this and the Leibniz rule it follows straightforwardly that for any $A, B, C \in \sy(N)$ there holds $\text{Cycle}\{\{A, B\}, C\} = 0$, so that $\{\,,\,\}$ is a Poisson structure on $\sy(N)$. By virtue of the Leibniz rule, the Poisson structure $\{\,,\,\}$ is completely determined by its action on vector fields. If $\nabla^{\prime}$ is a second torsion-free affine-connection, and $\{\,,\,\}^{\prime}$ is the Poisson structure determined by $\nabla^{\prime}$ according to \eqref{schouten}, then it is easy to check that for vector fields $A$ and $B$, $\{A, B\}^{\prime} = \{A, B\}$, from which it follows that the Poisson structures are the same, and therefore independent of the choice of torsion-free affine connection. Using this independence in a local coordinate chart it is easy to check that $\{\,,\,\}$ is simply the restriction of the tautological Poisson bracket on $C^{\infty}(T^{\ast}N)$. The preceeding shows that \eqref{schouten} gives a coordinate free way of computing this bracket, and this formulation will be used repeatedly in what follows. 

Each projective structure, $\en$, determines a projectively invariant, skew-symmetric, differential, bilinear pairing $\ppois: \emfb^{(i_{1}\dots i_{k_{1}})}[\lambda_{1}] \times \emfb^{(i_{1}\dots i_{k_{2}})}[\lambda_{2}] \to \emfb^{(i_{1}\dots i_{k_{1}+k_{2} - 1})}[\lambda_{1} + \lambda_{2}]$ defined by
\begin{align*}
\widetilde{\la a, b \ra}^{i_{1}\dots i_{k+l-1}} = \{\tilde{a}, \tilde{b}\}^{i_{1}\dots i_{k+l-1}}.
\end{align*}
This means that $\left(\widetilde{\la a, b\ra} - \{\tilde{a}, \tilde{b}\}\right)\wedge \eul = 0$. 
\begin{example}
If $\lambda_{1}, \lambda_{2}$ are non-excluded weights, for $a \in \emfb^{(i_{1}\dots i_{k_{1}})}[\lambda_{1}]$ and $b^{i} \in \emfb^{i}[\lambda_{2}]$, there hold
\begin{align}\label{pschouten1}
& \la a, b\ra^{i} = a^{p}\nabla_{p}b^{i} - b^{p}\nabla_{p}a^{i} + \tfrac{\lambda_{1}}{(n+1+\lambda_{2})}a^{i}\nabla_{p}b^{p} - \tfrac{\lambda_{2}}{(n+1+\lambda_{1})}b^{i}\nabla_{p}a^{p}\\
&\notag \la a, b\ra^{ij} = 2a^{p(i}\nabla_{p}b^{j)} - b^{p}\nabla_{p}a^{ij} + \tfrac{\lambda_{1}}{(\lambda_{2} + n + 1)}a^{ij}\nabla_{p}b^{p} - \tfrac{2\lambda_{2}}{(\lambda_{1} + n + 3)}b^{(i}\nabla_{p}a^{j)p}.
\end{align}
If $\lambda_{1} = 0 = \lambda_{2}$ then \eqref{pschouten1} shows that on vector fields the pairing $\ppois$ is the usual Schouten bracket, and so in this case $\ppois = \pois$, and in this case the notation $\pois$ will be used instead of $\ppois$.
\end{example}
\begin{lemma}
If $\lambda_{\ga}$ are non-excluded weights and $a_{\ga} \in \emfb^{(i_{1}\dots i_{k_{i}})}[\lambda_{\ga}]$ then
\begin{align*}
\la a_{1}\odot a_{2}, a_{3}\ra - a_{1}\odot \la a_{2}, a_{3}\ra - a_{2} \odot \la a_{1}, a_{3}\ra = \tfrac{\lambda_{3}(k_{1} + k_{2})}{(\lambda_{1} + \lambda_{2} + n + 2(k_{1}+k_{2}) - 1)}a_{3} \odot \op_{k_{1}+k_{2} - 1}(a_{1}, a_{2}).
\end{align*}
\end{lemma}
\begin{proof}
By Proposition \ref{liftsymprop} and the fact that $\pois$ is a Poisson bracket on $T^{\ast}\form$,
\begin{align*}
&\{\widetilde{a_{1}\odot a_{2}}, \tilde{a}_{3}\} -  \tilde{a}_{1} \odot \{\tilde{a}_{2}, \tilde{a}_{3}\} - \tilde{a}_{2}\odot\{\tilde{a}_{1}, \tilde{a}_{3}\} \\
&= \sum_{s = 1}^{k_{1}+k_{2}} \frac{\binom{k_{1}+k_{2}}{s}}{(\lambda_{1} + \lambda_{2} + n + 2(k_{1} + k_{2})-s)_{(s)}}\{\eul \odot \dots \odot \eul \odot \widetilde{\op_{k_{1}+k_{2}-s}(a_{1}, a_{2})}, \tilde{a}_{3}\}.
\end{align*}
To prove the claim it suffices to examine the $s = 1$ term. Because $\pois$ is a Poisson bracket,
\begin{align*}
\{\eul \odot \widetilde{\op_{k_{1} + k_{2} - 1}(a_{1}, a_{2})}, \tilde{a}_{3}\} = \eul \odot \{ \widetilde{\op_{k_{1} + k_{2} - 1}(a_{1}, a_{2})}, \tilde{a}_{3}\} +  \widetilde{\op_{k_{1} + k_{2} - 1}(a_{1}, a_{2})} \odot \{\eul, \tilde{a}_{3}\}
\end{align*}
and the term $\{\eul, \tilde{a}_{3}\} = \lambda_{3}\tilde{a}_{3}$, and from this the claim follows.
\end{proof}

The notation $O(m, f)$ (resp. $O(m, \tilde{f})$) is a shorthand for the phrase `terms of order at most $m$ in $f$ (resp. $\tilde{f}$)'

\begin{proposition}\label{skewstar}
If $\lambda_{1}$ and $\lambda_{2}$ are not excluded weights and $a_{p}^{i_{1}\dots i_{k_{p}}} \in \emfb^{(i_{1}\dots i_{k_{p}})}[\lambda_{p}]$ for $p = 1, 2$, then for any $\mu$, $[\op(a_{1}), \op(a_{2})] - \op(\la a_{1}, a_{2}\ra) \in \dop_{\mu, \mu+\Lambda}^{K-2}$.
\end{proposition}

\begin{proof}
By definition $\widetilde{\op_{0}(a_{1}, f)} = \tilde{a}_{1}^{I_{1}\dots I_{k}}\hnabla_{I_{1}\dots I_{k}}\tilde{f}$. Using the fact that $\hnabla_{[I}\hnabla_{J]}$ is a curvature term (so of $0$th order),
\begin{align*}
&\widetilde{\op_{0}(a_{2}, \op_{0}(a_{1}, f))} = \tilde{a}_{2}^{I_{k_{1}+1}\dots I_{K}}\tilde{a}_{1}^{I_{1}\dots I_{k_{1}}}\hnabla_{I_{k_{1}+1}}\dots\hnabla_{I_{K}}\hnabla_{I_{1}}\dots\hnabla_{I_{k_{1}}}\tilde{f}\\
& + k_{2}\tilde{a_{2}}^{II_{k_{1}+1}\dots I_{K-1}}\hnabla_{I}\tilde{a}_{1}^{I_{1}\dots I_{k_{1}}}\hnabla_{I_{k_{1}+1}}\dots\hnabla_{I_{K-1}}\hnabla_{I_{1}}\dots\hnabla_{I_{k_{1}}}\tilde{f} + O(K-2, \tilde{f})\\
& = \tilde{a}_{1}^{(I_{1}\dots I_{k_{1}}}\tilde{a}_{2}^{I_{k_{1}+1}\dots I_{K})}\hnabla_{I_{1}}\dots \hnabla_{I_{K}}\tilde{f} \\
& + k_{2}\tilde{a}_{2}^{I(I_{k_{1}+1}\dots I_{K-1}}\hnabla_{I}\tilde{a}_{1}^{I_{1}\dots I_{k_{1}})}\hnabla_{I_{1}}\dots\hnabla_{I_{K-1}}\tilde{f} + O(K-2, \tilde{f}),
\end{align*}
and from this there follows
\begin{align*}
&\widetilde{[\op(a_{1}), \op(a_{2})](f)} = \{\tilde{a}_{1}, \tilde{a}_{2}\}^{I_{1}\dots I_{K-1}}\hnabla_{I_{1}}\dots \hnabla_{I_{K-1}}\tilde{f} +  O(K-2, \tilde{f})\\
& = \widetilde{\la a_{1}, a_{2}\ra}^{I_{1}\dots I_{K-1}}\hnabla_{I_{1}}\dots\hnabla_{ I_{K-1}}\tilde{f} +O(K-2, \tilde{f})
= \widetilde{\op(\la a_{1}, a_{2}\ra)(f)} + O(K-2, \tilde{f}),
\end{align*}
and from this the claim follows.
\end{proof}

Formally extend the quantization map, $\op:\alg_{\lambda} \to \dop_{\lambda, \mu}$, to the algebra $\alg[[\ih]]$ by defining $\hop(a) = \ih^{k}\op(a)$ for $a \in \emfb^{(i_{1}\dots i_{k})}[\lambda]$. For weights $\lambda$ and $\mu$ which are not excluded, define $\hstar_{\lambda, \sigma}^{\mu}:\alg_{\lambda}\times \alg_{\sigma} \to \alg_{\lambda + \sigma}$ by $\hop(A \hstar_{\lambda,\sigma}^{\mu} B) = \hop(A)\circ \hop(B)$ and extend by linearity to $\alg_{\lambda}[[\ih]]$. By definition the product $\hstar_{\lambda, \sigma}^{\mu}$ is associative in the sense that 
\begin{align*}
a\hstar_{\lambda_{1}, \lambda_{2}+\lambda_{3}}^{\mu}(b \hstar_{\lambda_{2}, \lambda_{3}}^{\mu} c) = (a \hstar_{\lambda_{1}, \lambda_{2}}^{\mu} c)\hstar_{\lambda_{1}+ \lambda_{2}, \lambda_{3}}^{\mu} c.
\end{align*}
In particular, when $\lambda = \sigma = 0$, this product makes $\alg[[\ih]]$ into an associative algebra; in this case the product will be denoted by $\hstar^{\mu}$. In the particularly interesting case $\mu = -(n+1)/2$, $\hstar^{-(n+1)/2}$ will be written simply $\hstar$. Specializing $\ih = 1$ gives an associative multiplication, $\star$, on $\alg$. Denote the symmetric product of operators by $\odot$, so that $\op(a)\odot \op(b) = \tfrac{1}{2}(\op(a)\circ\op(b) + \op(b) \circ \op(a))$.

\begin{theorem}\label{starproducttheorem}
For $a \in \emfb^{(i_{1}\dots i_{k})}$, $b \in \emfb^{(i_{1}\dots i_{l})}$, and any $\mu \in \rea$,
\begin{align}\label{mainequation}
\op(a) \odot \op(b) - \op(a \odot b) - \tfrac{(n+1+2\mu)(k+l)}{2(n + 2(k+l) - 1)}\op(\op_{k+l-1}(a, b)) \in \dop_{\mu, 0}^{k+l-2},
\end{align}
and, as a consequence,
\begin{align}
&\label{ahmub} a \hmustar b = a \odot b + \tfrac{\ih}{2}\left (\{a, b\} +  \tfrac{(n+1+2\mu)(k+l)}{(n + 2(k+l) - 1)}\op_{k+l-1}(a, b)\right) + O(\ih^{2}),\\
&\label{liestar} a \hmustar b - b \hmustar a = \ih \{a, b\} + O(\ih^{2}).
\end{align}
In particular, when $\mu = -\tfrac{(n+1)}{2}$, 
\begin{align}
&\label{starproductproof}\op(a) \odot \op(b) - \op(a \odot b) \in \dop_{-(n+1)/2, 0}^{k+l-2},& &a\hstar b = a\odot b + \tfrac{\ih}{2}\{a, b\} + O(\ih^{2}).
\end{align}
\end{theorem}

\begin{proof}
To begin with assume $a \in \emfb^{(i_{1}\dots i_{k})}[\lambda]$, $b \in \emfb^{(i_{1}\dots i_{l})}[\sigma]$, $f \in \emfb[\mu]$, with $\lambda$, $\sigma$, and $\mu$ arbitrary (non-excluded) weights. The computations are involved, and it is easier to avoid mistakes if one does not specialize until the end to the particular weights appearing in the statement of the theorem.

The following identity will be used constantly, often without comment:
\begin{align}\label{dropterms1}
C^{I_{1}\dots I_{s}}\hnabla_{I_{1}}\dots \hnabla_{I_{s}}\tilde{f} = C^{(I_{1}\dots I_{s})}\hnabla_{I_{1}}\dots \hnabla_{I_{s}}\tilde{f} + O(s-2, \tilde{f}).
\end{align}
This follows from the fact that $\hnabla_{[I}\hnabla_{J]}$ is zeroth order. Likewise there follows from the proof of Proposition \ref{droptwoprop} the following, to be used frequently: 
\begin{align}\label{dropterms2}
\hnabla_{i_{1}}\dots \hnabla_{i_{s}}\tilde{f} = \widetilde{\nabla_{i_{1}}\dots \nabla_{i_{s}}f} + O(s-2, f).
\end{align}
By definition, \eqref{dropterms1}, \eqref{dropterms2}, and \eqref{iinfmix},
\begin{align*}
  & \widetilde{\op_{0}(a\odot b, f)} = \widetilde{a\odot b}^{I_{1}\dots I_{k+l}}\hnabla_{I_{1}}\dots \hnabla_{I_{k+l}}\tilde{f}\\
& = \widetilde{a\odot b}^{i_{1}\dots i_{k+l}}\widetilde{\nabla_{i_{1}}\dots \nabla_{i_{k+l}}f} \\
&+ (k+l)(\mu - (k+l-1))\widetilde{a\odot b}^{\infty i_{1}\dots i_{k+l-1}}\widetilde{\nabla_{i_{1}}\dots \nabla_{i_{k+l-1}}f} + O(k+l-2, f)
\end{align*}
By definition of the invariant lift,
\begin{align*}
\widetilde{a\odot b}^{\infty i_{1}\dots i_{k+l-1}} = -\tfrac{1}{(\lambda + \sigma + n + 2(k+l) - 1)}\widetilde{\nabla_{p}(a\odot b)^{i_{1}\dots i_{k+l-1}p}},
\end{align*}
so
\begin{align*}
& \op_{0}(a\odot b, f)= (a\odot b)^{i_{1}\dots i_{k+l}}\nabla_{i_{1}}\dots \nabla_{i_{k+l}}f \\
&+ \tfrac{(k+l)(k+l-1 - \mu)}{(\lambda + \sigma + n + 2(k+l)-1)}(\nabla_{p}(a\odot b)^{i_{1}\dots i_{k+l-1}p})\nabla_{i_{1}}\dots \nabla_{i_{k+l-1}}f + O(k+l-2, f)
\end{align*}
Now note that
\begin{align*}
a^{(i_{1}\dots i_{k}}\nabla_{p}b^{i_{k+1}\dots i_{k+l-1}p)} = \tfrac{l}{k+l}a^{(i_{1}\dots i_{k}}\nabla_{p}b^{i_{k+1}\dots i_{k+l-1})p} + \tfrac{k}{k+l}a^{p(i_{1}\dots i_{k-1}}\nabla_{p}b^{i_{k}\dots i_{k+l})}, 
\end{align*}
so that
\begin{align*}
&\nabla_{p}(a\odot b)^{i_{1}\dots i_{k+l-1}p} = a^{(i_{1}\dots i_{k}}\nabla_{p}b^{i_{k+1}\dots i_{k+l-1}p)} + b^{(i_{1}\dots i_{l}}\nabla_{p}a^{i_{l+1}\dots i_{k+l-1}p)}\\
& = \tfrac{l}{k+l}a^{(i_{1}\dots i_{k}}\nabla_{p}b^{i_{k+1}\dots i_{k+l-1})p} + \tfrac{k}{k+l}a^{p(i_{1}\dots i_{k-1}}\nabla_{p}b^{i_{k}\dots i_{k+l-1})}\\
& + \tfrac{k}{k+l}b^{(i_{1}\dots i_{l}}\nabla_{p}a^{i_{l+1}\dots i_{k+l-1})p} + \tfrac{l}{k+l}b^{p(i_{1}\dots i_{l-1}}\nabla_{p}a^{i_{l}\dots i_{k+l-1})},
\end{align*}
from which follows
\begin{align}
&\label{ingred2} \op_{0}(a\odot b, f)= (a\odot b)^{i_{1}\dots i_{k+l}}\nabla_{i_{1}}\dots \nabla_{i_{k+l}}f \\
&\notag + \tfrac{l(k+l-1 - \mu)}{(\lambda + \sigma + n + 2(k+l)-1)}\left(a^{(i_{1}\dots i_{k}}\nabla_{p}b^{i_{k+1}\dots i_{k+l-1})p} + b^{p(i_{1}\dots i_{l-1}}\nabla_{p}a^{i_{l}\dots i_{k+l-1})} \right) \nabla_{i_{1}}\dots \nabla_{i_{k+l-1}}f \\
&\notag + \tfrac{k(k+l-1 - \mu)}{(\lambda + \sigma + n + 2(k+l)-1)}\left(a^{p(i_{1}\dots i_{k-1}}\nabla_{p}b^{i_{k}\dots i_{k+l-1})}
 + b^{(i_{1}\dots i_{l}}\nabla_{p}a^{i_{l+1}\dots i_{k+l-1})p} \right) \nabla_{i_{1}}\dots \nabla_{i_{k+l-1}}f \\ &\notag + O(k+l-2, f)
\end{align}
Note that the specialization $\lambda = 0 = \sigma$ and $\mu = -(n+1)/2$ gives $\tfrac{(k+l-1 - \mu)}{(\lambda + \sigma + n + 2(k+l)-1)} = 1/2$. This completes the computation of the top order terms of $\op_{0}(a\odot b, f)$. Now consider the top order terms of $\op_{0}(a, \op_{0}(b, f))$. Using \eqref{dropterms1} yields
\begin{align*}
&\widetilde{\op_{0}(a, \op_{0}(b, f))} = \tilde{a}^{(I_{1}\dots I_{k}}\tilde{b}^{I_{k+1}\dots I_{k+l})}\hnabla_{I_{1}}\dots \hnabla_{I_{k+l}}\tilde{f} \\
& + k\tilde{a}^{QI_{1}\dots I_{k-1}}(\hnabla_{Q}\tilde{b}^{I_{k}\dots I_{k+l-1}})\hnabla_{I_{1}}\dots\hnabla_{I_{k+l-1}}\tilde{f} + O(k+l -2, \tilde{f}).
\end{align*}
Observe that 
\begin{align*}
&\tilde{a}^{(\infty i_{1}\dots i_{k-1}}\tilde{b}^{i_{k}\dots i_{k+l-1})} = \tfrac{k}{k+l}\tilde{a}^{\infty (i_{1}\dots i_{k-1}}\tilde{b}^{i_{k}\dots i_{k+l-1})} + \tfrac{l}{k+l}\tilde{b}^{\infty (i_{1}\dots i_{l-1}}\tilde{a}^{i_{l}\dots i_{k+l-1})}\\
& = -\tfrac{k}{(k+l)(\lambda + n + 2k-1)}\widetilde{b^{(i_{1}\dots i_{l}}\nabla_{p}a^{i_{l+1}\dots i_{k+l-1})p}}  -\tfrac{l}{(k+l)(\sigma + n + 2l-1)}\widetilde{a^{(i_{1}\dots i_{k}}\nabla_{p}b^{i_{k+1}\dots i_{k+l-1})p}}. 
\end{align*}
With \eqref{dropterms2} and \eqref{iinfmix}, this yields
\begin{align}
&\label{step4}\tilde{a}^{(I_{1}\dots I_{k}}\tilde{b}^{I_{k+1}\dots I_{k+l})}\hnabla_{I_{1}}\dots \hnabla_{I_{k+l}}\tilde{f} = \widetilde{a^{(i_{1}\dots i_{k}}b^{i_{k+1}\dots i_{k+l})}}\widetilde{\nabla_{i_{1}}\dots \nabla_{i_{k+l}}f}\\
&\notag + \tfrac{k(k+l-1-\mu)}{(\lambda + n + 2k -1)}\widetilde{b^{(i_{1}\dots i_{l}}\nabla_{p}a^{i_{l+1}\dots i_{k+l-1})p}}\\
&\notag +  \tfrac{l(k+l-1-\mu)}{(\sigma + n + 2l -1)}\widetilde{a^{(i_{1}\dots i_{k}}\nabla_{p}b^{i_{k+1}\dots i_{k+l-1})p}}  + O(k+l-2, f).
\end{align}
Next are analyzed the contributions of the term $ k\tilde{a}^{QI_{1}\dots I_{k-1}}(\hnabla_{Q}\tilde{b}^{I_{k}\dots I_{k+l-1}})\hnabla_{I_{1}}\dots\hnabla_{I_{k+l-1}}\tilde{f}$, which equals
\begin{align*}
k\left(\tilde{a}^{qi_{1}\dots i_{k-1}}\hnabla_{q}\tilde{b}^{i_{k}\dots i_{k+l-1}} + \tilde{a}^{\infty i_{1}\dots i_{k-1}}\hnabla_{\infty}\tilde{b}^{i_{k}\dots i_{k+l-1}}\right)\widetilde{\nabla_{i_{1}}\dots \nabla_{i_{k+l-1}}f} + O(k+l-2, f).
\end{align*}
Using \eqref{intermed1} gives
\begin{align}
\label{fs1}& (\tilde{a}^{qi_{1}\dots i_{k-1}}\hnabla_{q}\tilde{b}^{i_{k}\dots i_{k+l-1}})\widetilde{\nabla_{i_{1}}\dots \nabla_{i_{k+l-1}}f} =&\\
\notag  ( a^{p(i_{1}\dots i_{k-1}}\nabla_{b}b^{i_{k}\dots i_{k+l-1})} &- \tfrac{l}{(\sigma + n + 2l - 1)}a^{(i_{1}\dots i_{k}}\nabla_{p}b^{i_{k+1}\dots i_{k+l-1})p})\widetilde{\nabla_{i_{1}}\dots \nabla_{i_{k+l-1}}f} + O(k+l-2, f). &
\end{align}
Using \eqref{induct1} and \eqref{intermed0} gives 
\begin{align*}
&k\tilde{a}^{\infty i_{1}\dots i_{k-1}}\hnabla_{\infty}\tilde{b}^{i_{k}\dots i_{k+l-1}} = -\tfrac{k(l + \sigma)}{(\lambda + n + 2k - 1)}b^{i_{1}\dots i_{l}}\nabla_{p}a^{i_{l+1}\dots i_{k+l-1}p},\end{align*}
so that
\begin{align}
&\label{fs2} k\tilde{a}^{\infty i_{1}\dots i_{k-1}}\hnabla_{\infty}\tilde{b}^{i_{k}\dots i_{k+l-1}}\widetilde{\nabla_{i_{1}}\dots \nabla_{i_{k+l-1}}f}\\
&\notag  = -\tfrac{k(l + \sigma)}{(\lambda + n + 2k - 1)}\widetilde{b^{(i_{1}\dots i_{l}}\nabla_{p}a^{i_{l+1}\dots i_{k+l-1})p}}\widetilde{\nabla_{i_{1}}\dots \nabla_{i_{k+l-1}}f} + O(k+l-2, f).
\end{align}
Adding \eqref{step4}, \eqref{fs1}, and \eqref{fs2} yields
\begin{align*}
&\op_{0}(a, \op_{0}(b, f)) = a^{(i_{1}\dots i_{k}}b^{i_{k+1}\dots i_{k+l})}\nabla_{i_{1}}\dots \nabla_{i_{k+l}}f\\
& +\tfrac{k(k-1-\mu-\sigma)}{(\lambda + n + 2k -1)}b^{(i_{1}\dots i_{l}}\nabla_{p}a^{i_{l+1}\dots i_{k+l-1})p}\nabla_{i_{1}}\dots \nabla_{i_{k+l-1}}f\\
& + \tfrac{l(l-1-\mu)}{(\sigma + n + 2l -1)}a^{(i_{1}\dots i_{k}}\nabla_{p}b^{i_{k+1}\dots i_{k+l-1})p}\nabla_{i_{1}}\dots \nabla_{i_{k+l-1}}f\\
& + ka^{p(i_{1}\dots i_{k-1}}\nabla_{p}b^{i_{k}\dots i_{k+l-1})} \nabla_{i_{1}}\dots \nabla_{i_{k+l-1}}f + O(k-l-2, f).
\end{align*}
Now symmetrizing yields
\begin{align}
&\label{ingred3} \op(a)\odot\op(b)(f) -  a^{(i_{1}\dots i_{k}}b^{i_{k+1}\dots i_{k+l})}\nabla_{i_{1}}\dots \nabla_{i_{k+l}}f\\
&\notag = \left(\tfrac{l}{2}\tfrac{(2l -2 - 2\mu - \lambda)}{\sigma + n + 2k -1)}a^{(i_{1}\dots i_{k}}\nabla_{p}b^{i_{k+1}\dots i_{k+l-1})p} + \tfrac{l}{2}b^{p(i_{1}\dots i_{l-1}}\nabla_{p}a^{i_{l}\dots i_{k+l-1})} \right)\nabla_{i_{1}}\dots\nabla_{i_{k+l-1}}f\\
&\notag + \left(\tfrac{k}{2}\tfrac{(2k-2-2\mu - \sigma)}{(\lambda + n + 2k -1)}b^{(i_{1}\dots i_{l}}\nabla_{p}a^{i_{l+1}\dots i_{k+l-1})p} + \tfrac{k}{2}a^{p(i_{1}\dots i_{k-1}}\nabla_{p}b^{i_{k}\dots i_{k+l-1})} \right)\nabla_{i_{1}}\dots\nabla_{i_{k+l-1}}f\\
&\notag + O(k+l-2, f).
\end{align}
Together with the specialization $\lambda = 0 = \sigma$, \eqref{opkminusone}, \eqref{ingred2}, and \eqref{ingred3} give \eqref{mainequation}. With Proposition \ref{skewstar}, \eqref{mainequation} implies \eqref{ahmub} and \eqref{liestar}. 
\end{proof}

\begin{corollary}[Corollary of Theorem \ref{starproducttheorem}]
The multiplication $\hmustar$ is a graded differential star product on $\alg$ adapted to the tautological Poisson structure on $\alg$ and having the form \eqref{speq}.
\end{corollary}

\begin{proof}
By Theorem \ref{starproducttheorem} and Proposition \ref{skewstar}, 
\begin{align*}
\dop = \op(a)\circ \op(b) - \op(a \odot b) - \tfrac{1}{2}\op(\{ a, b\}) - \tfrac{(k+l)(n+1+2\mu)}{2(n + 2(k+l)-1)}\op_{k+l-1}(a, b)
\end{align*}
is a differential operator of order at most $k + l - 2$. Define $B_{2}(a, b) \in \emfb^{(i_{1}\dots i_{k+l-2})}$ to be its principal symbol. Then $\dop - \op(B_{2}(a, b))$ is a differential operator of order at most $k + l - 3$, and let $B_{3}(a, b) \in \emfb^{(i_{1}\dots i_{k+l-3})}$ be its principal symbol. Proceeeding in this way shows that there are $B_{r}(a, b) \in \emfb^{(i_{1}\dots i_{k+l-r})}$ such that $\op(a)\circ \op(b) - \op(a \odot b) = \sum_{r = 1}^{k+l}\op(B_{r}(a, b))$. The tensors $B_{r}(a, b)$ are evidently bilinear differential operators in $a$ and $b$. Since the differential operators involved in defining the $B_{r}$ are all projectively invariant, the operators $B_{r}$ will be also. From the definition of $\op(a) \circ \op(b)$ it is straightforward to see that $B_{r}(a, b)$ has order at most $r$ in each of $a$ and $b$. That $\hmustar$ satisfies \eqref{starproductdefined} now follows by a formal calculation. That the constant function $1$ is a unit for $\hmustar$ is immediate. By \eqref{liestar} the Lie bracket of $\hmustar$ agrees with the usual Poisson structure on $\alg$ to first order.

By \eqref{ahmub} of Theorem \ref{starproducttheorem} the Hochshild $2$-cochain of $\alg$ with coefficients in $\alg$ defined on $\alg_{k} \times \alg_{l}$ by $H(a, b) = \tfrac{k+l}{n + 2(k+l)- 1}\op_{k+l-1}(a, b)$ must be a coboundary. If a Hochschild $1$-cochain, $C_{1}$, is defined on $\alg_{k}$ by $C_{1}(a) = \tfrac{k}{n+2k-1}\nabla_{p}a^{i_{1}\dots i_{k-1}p}$, then it is easy to check that $\partial C_{1}(a, b) = a \odot C_{1}(b) - C_{1}(a \odot b) + C_{1}(a) \odot b = H$. Note that choosing a different representative of $\en$ only adds to $C_{1}$ a $1$-cocycle of the form $a^{i_{1}\dots i_{k}} \to k\gamma_{p}a^{i_{1}\dots i_{k-1}p}$, so that $\partial C_{1}$ is uniquely defined. The cochain $C_{1}$ can be written as a formal pseudo-differential operator on $\alg$ as follows. Let $D$ be the divergence operator on $T^{\ast}M$ associated to the representative $\nabla \in \en$ as in Lemma \ref{divlemma}. Then $C_{1} = \tfrac{1}{n+1 + 2\lie_{\teul}}\circ D$, where $\teul$ is the Euler vector field on $T^{\ast}M$, and $\tfrac{1}{n+1+2\lie_{\teul}}$ is intepreted as a formal power series in $\lie_{\teul}$. Changing the choice of representative of $\en$ only perturbs $C_{1}$ by adding a $1$-cocycle, and so $\partial C_{1}$ is well-defined, and evidently a formal bilinear pseudo-differential operator. 
\end{proof}
It follows also that the star product $A \bhstar B = T(T^{-1}(A) \hmustar T^{-1}(B))$ gauge-equivalent to $\hmustar$ by the change of gauge $T = Id + \ih \tfrac{(n+1+2\mu)}{2}C_{1}$ satisfies $B_{1} = \tfrac{1}{2}\pois$.

\eqref{starvectora} and \eqref{starvectorb} suggest that in general $B^{\mu}_{r}(a, b)$ should be expressible in terms of sums of iterated expressions composed of $\op_{\be}$ and $\pois$.

\begin{proposition}\label{symmetricstarproduct}
Writing $\Bar{\ih} = -\ih$, there holds for $a \in \emfb^{(i_{1}\dots i_{k})}$ and $b \in \emfb^{(i_{1}\dots i_{l})}$,
\begin{align}\label{skewchmustar}
a\hmustar b = b \chmustar a.
\end{align}
In particular the star product $\hstar$ is symmetric.
\end{proposition}
\begin{proof}
Let $u \in \emfb_{c}[\mu]$ and $v \in \emfb_{c}[-\mu-n-1]$. Write $a \hmustar b = \sum_{r = 0}^{k+l}\ih^{r}B_{r}^{\mu}(a, b)$ and $a \hcmustar b = \sum_{r = 0}^{k+l}\ih^{r}B_{r}^{-\mu-n-1}(a, b)$. The claim amounts to showing $B_{r}^{\mu}(a, b) = (-1)^{r}B_{r}^{-\mu-n-1}(b, a)$, and so when $\mu = -(n+1)/2$ shows the symmetry of $\hstar$. Now $\hop(a \bhmustar b) = \ih^{k+l}\sum_{r = 0}^{k+l}(-1)^{r}\op(B_{r}^{\mu}(a, b))$, so using Theorem \ref{selfadjoint} it makes sense to write 
\begin{align*}
&(\hop(a \bhmustar b)(u), v) = \ih^{k+l}\sum_{r = 0}^{k+l}(-1)^{r}(\op(B_{r}^{\mu}(a, b))(u), v) = \ih^{k+l}\sum_{r = 0}^{k+l}(-1)^{k+l}(u, \op(B_{r}^{\mu}(a, b))(v))\\
& = (-1)^{k+l}\sum_{r = 0}^{k+l}\ih^{r}(u, \hop(B_{r}^{\mu}(a, b))(v)) = (-1)^{k+l}(u, \hop(a \hmustar b)(v)).
\end{align*}
On the other hand, by definition of $\hmustar$ and Theorem \ref{selfadjoint},
\begin{align*}
&(\hop(a \hmustar b)(u), v) = (\hop(a)\circ \hop(b)(u), v) = (-1)^{k}(\hop(b)(u), \hop(a)(v)) \\
&= (-1)^{k+l}(u, \hop(b) \circ \hop(a)(v)) = (-1)^{k+l}(u, \hop(b \chmustar a)(v)).
\end{align*}
Hence the equality $(u, \hop(a \hmustar b)(v)) = (u, \hop(b \chmustar a)(v))$ holds for all choices of $u$ and $v$, which implies the claim.
\end{proof}
In one dimension Proposition \ref{symmetricstarproduct} specializes to equation 4.6 of \cite{Cohen-Manin-Zagier}.
\begin{example}
Let $a^{i} \in \emfb^{i}[\lambda_{1}]$, $b^{i} \in \emfb^{i}[\lambda_{2}]$, and $f \in \emfb[\mu]$. Write $\Lambda = \lambda_{1}+\lambda_{2}$. Then
\begin{align}
\label{op2inter1} & \op(a)\odot \op(b) - \op(a\odot b) = \tfrac{(\Lambda + n + 2\mu +1))}{(\Lambda + n + 3)}\op(\op_{1}(a, b)) + \tfrac{\mu(\Lambda + \mu + n + 1)}{(\Lambda + n+2)_{(2)}}\op(\op_{0}(a, b)).
\end{align}
In particular, if $\Lambda = 0$ and $\mu = - (n+1)/2$, then
\begin{align*}
\op(a)\circ \op(b) = \op(a\odot b) + \tfrac{1}{2}\op(\{ a, b\}) - \tfrac{(n+1)}{4(n+2)}\op(\op_{0}(a, b)).
\end{align*}
This shows that 
\begin{align}
&\label{ahmustarb} a\hmustar b = a\odot b + \ih\left(\tfrac{1}{2}\{ a, b\} + \tfrac{n+1 + 2\mu}{n+3}\op_{1}(a, b)\right ) +  \ih^{2}\tfrac{\mu(\mu + n+1)}{(n+2)(n+1)}\op_{0}(a, b),\\
&\label{starvector} a\hstar b = a\odot b + \tfrac{\ih}{2}\{ a, b\} - (\tfrac{\ih}{2})^{2}\tfrac{(n+1)}{(n+2)}\op_{0}(a, b).
\end{align} 
To show \eqref{op2inter1} it suffices to verify the identities,
\begin{align*}
&\hnabla_{I}\tilde{a}^{(I}\tilde{b}^{J)} = \widetilde{\op_{1}(a, b)}^{J} + \frac{1}{\Lambda + n +1}\widetilde{\op_{0}(a, b)}\eul^{J},&
&\widetilde{\la a, b\ra}^{I} = \{\tilde{a}, \tilde{b}\}^{I},
\end{align*}
and to use the definitions of the quantities involved.
\end{example}

\begin{corollary}
The star product, $\hmustar$, is \textbf{covariant} with respect to $\diff(M)$ in the sense that, for $a^{i}, b^{i} \in \emfb^{i}$ there holds $a \hmustar b - b \hmustar a = \ih\{a, b\}$. That is, the restriction to vector fields of the star product commutator agrees with $\ih\pois$.
\end{corollary}

\begin{example}\label{starvectorexample}
Let $a^{i_{1}\dots i_{k}} \in \emfb^{(i_{1}\dots i_{k})}$ and $f \in \emfb[0]$. Computing $\op(fa)(u)$ using \eqref{liftsym} and \eqref{iinfmix} yields 
\begin{align}\label{starvector0}
f \hmustar a = fa + \sum_{r \geq 1}\ih^{r}\tfrac{\binom{k}{s}\binom{\mu -k + s}{s}}{\binom{n+2k-s}{s}}\op_{k-s}(a, f).
\end{align}
Using $\{a, f\} = k\op_{k-1}(a, f)$ in \eqref{starvector0} yields \eqref{starvectorb}, and using $a \hmustar f = f\,_{-\epsilon}{\star}^{-\mu-n-1}\,\,a$ with \eqref{starvectorb} yields \eqref{starvectora}.
\begin{align}
&\label{starvectorb} f \hmustar a = fa +\tfrac{\ih}{2}\left(\{f, a\} + \tfrac{k(n+1+2\mu)}{n+2k-1}\op_{k-1}(f, a)\right)+\sum_{s = 2}^{k}\ih^{s}\frac{\binom{k}{s}\binom{\mu-k+s}{s}}{\binom{n + 2k - s}{s}} \op_{k-s}(a,f),\\
&\label{starvectora} a\hmustar f = fa + \tfrac{\ih}{2}\left(\{a, f\} + \tfrac{k(n+1+2\mu)}{n+2k-1}\op_{k-1}(f, a)\right) + \sum_{s = 2}^{k}\ih^{s}\frac{\binom{k}{s}\binom{n+k+\mu}{s}}{\binom{n+2k-s}{s}} \op_{k-s}(a,f).
\end{align}
It follows that
\begin{align}
\label{starvectorc} \tfrac{1}{2}(a \hstar f - f \hstar a - \ih\{a, f\}) = \sum_{s = 1}^{\lfloor k/2 \rfloor}\ih^{2s+1}\frac{\binom{k}{2s+1}\binom{\frac{2k+n-1}{2}}{2s+1}}{\binom{n+2k - 2s - 1}{2s+1}}\op_{k -2s - 1}(a, f).
\end{align}
In the special case that $a^{ij} \in \emfb^{(ij)}$ and $f \in \emfb[0]$ there holds
\begin{align}
&\label{starvector2} a\hstar f = fa + \tfrac{\ih}{2}\{ a, f\} + (\tfrac{\ih}{2})^{2}\tfrac{(n+3)}{(n+2)}\op_{0}(a, f) = f \hminusstar a.
\end{align} 
The associativity of $\hstar$ and equations \eqref{starvector} and \eqref{starvector2} imply that for $a^{i}, b^{i} \in \emfb^{i}$ and$ f \in \emfb[0]$ there hold the identities
\begin{align*}
&\op_{0}(a, fb) - \op_{0}(b, fa) = -\tfrac{n+2}{n+1}\{\{a, b\}, f\}\\
&\op_{0}(a\odot b, f) = \tfrac{(n+1)}{2(n+3)}(2\op_{0}(a, b)f - \op_{0}(a, fb) - \op_{0}(b, fa)) + \tfrac{(n+2)}{2(n+3)}(\{\{f, a\}, b\} + \{\{f, b\}, a\}).
\end{align*}
It would be interesting to understand these identities as special cases of some general family of identities.

Applying \eqref{starvectora} in the one-dimensional case yields, for $u_{1} \in \emfb[2k]$ and $u_{2} \in \emfb[0]$,
\begin{align}\label{starvectorone}
u_{1} \hmustar u_{2} = u_{1}u_{2} + \sum_{s = 1}^{k}\ih^{s}\frac{(-1)^{s}\binom{k}{s}\binom{k+\mu + 1}{s}}{\binom{2k}{s}\binom{2k+1-s}{s}}\rc_{s}(u_{1}, u_{2}).
\end{align}
Let the notation be as in Section \ref{onedimsection}. In \cite{Cohen-Manin-Zagier} an associative multiplication, $\cmzmu(u_{1}, u_{2})$, of $u_{i} \in \emfb[2k_{i}]$, is defined by
\begin{align}
&\label{cmzmult}\cmzmu(u_{1}, u_{2}) = u_{1}u_{2} + \sum_{r \geq 1} t_{r}^{\mu}\rc_{r}(u_{1}, u_{2}),&\\
 &\notag t_{r}^{\mu} = (-1/4)^{r}\sum_{j \geq 0}\frac{\binom{r}{j}\binom{-1/2}{j}\binom{-3/2 -\mu}{j}\binom{1/2 + \mu}{j}}{\binom{k_{1} - 1/2}{j}\binom{k_{2}-1/2}{j}\binom{r - k_{1}-k_{2} - 3/2}{j}}.
\end{align}
The multiplications $\cmzmu$ and $\cmzmub$ are the same, so that $\cmzmu$ cannot be obtained from $\hmustar$ merely by specializing the parameter $\ih$. A multiplication such that $\hmubox = \hmubbox$ can be defined by $u_{1} \hmubox u_{2} = \tfrac{1}{2}(u_{1}\hmustar u_{2} + u_{1}\hmubstar u_{2})$. However, when $\mu \neq -1$, associativity is not clear. Using \eqref{starvectorone} it can be shown that in the special case $u_{1} \in \emfb[2k]$ and $u_{2} \in \emfb[0]$ that the coefficient of $\ih^{2}$ in $u_{1} \hmubox u_{2}$ equals $\tfrac{1}{4}t_{2}^{\mu}$, so $\hmubox$ has at least the right spirit. The remark at the very end of \cite{Duval-ElGradechi-Ovsienko} that for $\mu = -1$ the one-dimensional projectively invariant star products correspond to the multiplications obtained in \cite{Cohen-Manin-Zagier} is correct in the sense that, as just discussed, there must be an explicit relation between $\cmzmu$ and $\hmustar$. However, the form of such a relation does not appear to follow immediately from the uniqueness of homogeneous invariant star products proved in \cite{Duval-ElGradechi-Ovsienko}, because even for $\mu = -1$ it is not apparent that the multiplications $\cmzmu$ are star products. To pass from $\hmustar$ to $\cmzmu$ one has to perform manipulations involving specializing the parameter $\ih$, and the result of these manipulations is hidden in the explicit formula for $t_{r}^{\mu}$; if $t_{r}^{\mu}$ is replaced by $\ih^{r}t_{r}^{\mu}$ in \eqref{cmzmult} it is not clear that the resulting multiplication should be associative.
\end{example}

\subsection{Multiplication at Infinity}\label{infsection}
\begin{lemma}
For every $\mu \in \rea$, the restriction to $\alg_{k} \times \alg_{l}$ of $B_{r}^{\mu}$ is a polynomial in $\mu$ of degree at most $r$.
\end{lemma}
\begin{proof}
Let $a \in \alg_{k}$ and $f \in \emfb[\mu]$. It is claimed that $\op(a)(f)$ has the form $\sum_{s = 0}^{k}C_{s}^{i_{1}\dots i_{s}}\nabla_{i_{1}\dots i_{s}}f$, where $C_{s}$ is a universal polynomial in $a$ and its covariant derivatives which is polynomial in $\mu$ of degree at most $k-s$. By definition $\widetilde{\op(a)(f)} = \tilde{a}^{I_{1}\dots I_{k}}\hnabla_{I_{1}}\dots \hnabla_{I_{k}}\tilde{f}$. For $f \in \emfb[\mu]$, factors involving $\mu$ appear in this expression only through differentiation of $\tilde{f}$ in the vertical direction as in \eqref{iinfmix}. From \eqref{iinfmix} it is clear that a polynomial of degree $k -s$ in $\mu$ appears only as the coefficient of a term involving $s$ horizontal derivatives of $\tilde{f}$, and this shows the claim.

The proof of the lemma is by induction on $r$. The claim is true for $r = 1$ by \eqref{mainequation}. Let $a \in \alg_{k}$ and $b \in \alg_{l}$. By the inductive hypothesis, for $1 \leq s \leq r-1$, $B_{s}^{\mu}(a, b)$ is a polynomial in $\mu$ of degree at most $s$. $B_{r}^{\mu}$ is defined by the requirement that $B_{r}^{\mu}(a, b)$ be the principal symbol of the operator $\dop = \op(a)\circ \op(b) - \op(a \odot b) - \sum_{s = 1}^{r-1}\op(B_{s}^{\mu}(a, b))$, which has order at most $k+l - r$. By the inductive hypothesis and the remarks in the preceeding paragraph for $1 \leq s\leq r-1$ the coefficient of the term in $\op(B_{s}^{\mu}(a, b))(f)$ involving $k + l - r$ derivatives of $f$ is a polynomial in $\mu$ of degree at most $s + (r - s) = r$. Likewise, by the remarks in the preceeding paragraph the coefficient of the term in $\op(a) \circ \op(b)(f) - \op(a \odot b)(f)$ involving $k + l - r$ derivatives of $f$ is a polynomial in $\mu$ of degree at most $r$. $B_{r}^{\mu}(a, b)$ must agree with the coefficient in $\dop$ of $\nabla_{(i_{1}}\dots \nabla_{i_{k+l-r})}f$, which was just shown to be a polynomial in $\mu$ of degree at most $r$. This proves the lemma.
\end{proof}

\begin{proposition}
For $t\in \rea$, the limit of $\lim_{t\to 0} \,_{t}{\star}^{c/t}\,\, = \cinfstar$ exists and is a commutative, graded differential star product having the form \eqref{speq2}.
\end{proposition}
\begin{proof}
Since the restriction to $\alg_{k} \times \alg_{l}$ of $B_{r}^{\mu}$ is a polynomial in $\mu$ of degree at most $r$, for $t\in \rea$, $a \in \alg_{k}$ and $b \in \alg_{l}$, the limit $\lim_{t \to 0} a \,_{t}{\star}^{c/t}\,\, b$ exists and differs from $a \odot b$ by a polynomial in $c$ of order at most $k+l$. For each $c \in \rea$ define $a \cinfstar b = \lim_{t \to 0} a \,_{t}{\star}^{c/t}\,\, b$. Taking limits in \eqref{skewchmustar} shows that $a \cinfstar b = b \cinfstar a$, so that $\cinfstar$ defines a commutative, associative multiplication on $\alg$. For $a \in \alg_{k}$ and $b \in \alg_{l}$, the product $a \cinfstar b$ differs from the usual multiplication by a polynomial of degree at most $k+l$ in $c$, the coefficient of $c^{s}$ of which is a projectively invariant bilinear differential operator in $a$ and $b$ of degree at most $k+l-s$ in each argument. That $\cinfstar$ has the form \eqref{speq2} follows from \eqref{speq}.
\end{proof}
The limiting definition of $\cinfstar$ was motivated by the multiplication $\cmzinf$ defined in the one-dimensional case in section 4 of \cite{Cohen-Manin-Zagier}. This is explained further in Example \ref{infmultex} below.

\begin{example}
 For $a^{i}, b^{i} \in \emfb^{i}$, using \eqref{ahmustarb} yields
\begin{align*}
a \cinfstar b = a \odot b + \tfrac{2c}{n+3}\op_{1}(a, b) + \tfrac{c^{2}}{(n+2)(n+1)}\op_{0}(a, b).
\end{align*}
Comparing this example with \eqref{liftsymprop} suggests Theorem \ref{inftheorem}.
\end{example}

\begin{proof}[Proof of Theorem \ref{inftheorem}]
Let $a_{i} \in \emfb^{(i_{1}\dots i_{k_{i}})}$ and $f \in \emfb[\mu]$. From \eqref{liftsym} there follows
\begin{align}\label{infintermed1}
\tilde{a}_{1}^{(I_{1}\dots I_{k_{1}}}\tilde{a}_{2}^{I_{k_{1}+1}\dots I_{K})}\hnabla_{I_{1}}\dots \hnabla_{I_{K}}\tilde{f} = \widetilde{\op(a_{1}\odot a_{2})(f)} + \sum_{s = 1}^{K}\frac{\binom{K}{s}\binom{\mu + s - K}{s}}{\binom{n + 2K - s}{s}}\widetilde{\op(\op_{K-s}(a_{1}, a_{2}))(f)}
\end{align}
It is claimed that
\begin{align*}
\widetilde{\op(a_{1})\circ \op(a_{2})(f)} - \tilde{a}_{1}^{(I_{1}\dots I_{k_{1}}}\tilde{a}_{2}^{I_{k_{1}+1}\dots I_{K})}\hnabla_{I_{1}}\dots \hnabla_{I_{K}}\tilde{f}
\end{align*}
is a polynomial in $\mu$ of degree at most $K - 1$. If there occurs in either $\widetilde{\op(a_{1})\circ \op(a_{2})(f)}$ or $\tilde{a}_{1}^{(I_{1}\dots I_{k_{1}}}\tilde{a}_{2}^{I_{k_{1}+1}\dots I_{K})}\hnabla_{I_{1}}\dots \hnabla_{I_{K}}\tilde{f}$ a term of order $K$ in $\mu$, this term must involve $K$ vertical derivatives of $\tilde{f}$. With respect to any choice of projective scale the term in either expression involving $K$ vertical derivatives of $\tilde{f}$ is $\tilde{a}_{1}^{(\infty\dots \infty}\tilde{a}_{2}^{\infty\dots \infty)}\hnabla_{\infty}\dots \hnabla_{\infty}\tilde{f}$, and this shows the claim. With \eqref{infintermed1} and \eqref{speq} this claim shows that 
\begin{align}\label{infintermed2}
\op(a_{1})\circ \op(a_{2})(f) - \op(a_{1}\odot a_{2})(f) -  \sum_{s = 1}^{K}\frac{\binom{K}{s}\binom{\mu + s - K}{s}}{\binom{n + 2K - s}{s}}\op(\op_{K-s}(a_{1}, a_{2}))(f)
\end{align}
is a polynomial in $\mu$ of degree at most $K - 1$ and a differential operator in $f$ of degree at most $K - 1$. Setting $\mu = \frac{c}{t}$ and multiplying \eqref{infintermed2} by $t^{K}$ yields
\begin{align*}
&\thop(a_{1} \tstar a_{2}) = \thop(a_{1}\odot a_{2}) +  \sum_{s = 1}^{K}\frac{\binom{K}{s}\binom{c/t + s - K}{s}}{\binom{n + 2K - s}{s}}t^{s}\thop(\op_{K-s}(a_{1}, a_{2})) + O(t)\thop(b)\\
& = \thop(a_{1}\odot a_{2}) +  \sum_{s = 1}^{K}\frac{\binom{K}{s}}{(n + 2K - s)_{(s)}}c^{s}\thop(\op_{K-s}(a_{1}, a_{2})) + O(t)\thop(b^{\prime})
\end{align*}
where $b$ and $b^{\prime}$ have degree at most $K - 1$. This implies
\begin{align*}
a_{1} \tstar a_{2} = a_{1}\odot a_{2} +  \sum_{s = 1}^{K}\frac{\binom{K}{s}}{(n + 2K - s)_{(s)}}c^{s}\op_{K-s}(a_{1}, a_{2}) + O(t)
\end{align*}
and letting $t \to 0$ gives \eqref{binfr}. With this definition the associativity of $\cinfstar$ is not entirely clear. Associativity follows immediately from the following alternative defintion of $\cinfstar$, which shows that $\cinfstar$ is implicit already in \eqref{liftsym}.

The symmetric product of symmetric polyvectors on $\form$ having homogeneity $0$ is again homogeneous of degree $0$, so there is a graded subalgebra $\pol_{0}(T^{\ast}\form) \subset \pol(T^{\ast}\form)$. There is a distinguisshed element, $\heul \in \pol_{0}(T^{\ast}\form)$ determined by $\eul$ which in the notation of Section \ref{anothersection} can be written $\heul = \hPhi(\eul)$. An element of $\pol_{0}(T^{\ast}\form)$ is necessarily polynomial in $\heul$. The algebra, $\alg[c]$, of polynomials in $c$ with coefficients in $\alg$, is graded by the sum of the degree of the coefficients and the degree in $c$. The invariant lift of elements of $\alg$ determined by the projective structure $\en$ determines a graded linear isomorphism $\tlop:\alg[c] \to \pol_{0}(T^{\ast}\form)$ defined by letting $\tlop(ac^{l}) = \hPhi(\tilde{a})\heul^{l}$ for $a \in S^{k}(TM)$ and extending linearly. The map $\tlop$ is clearly injective. To see that it is surjective, let $A \in \pol_{0}(T^{\ast}\form)$ be homogeneous of degree $k$ in the fibers (i.e. correspond to a section of $S^{k}(T\form)$); then $A = \sum_{m = 0}^{k}\tlop(a_{m})\heul^{m} = \tlop(\sum_{m=0}^{k}a_{m}c^{m})$ where the $a_{i} \in S^{i}(TM)$ are uniquely determined inductively by the requirement that for $0 \leq i \leq k-1$ there hold $A - \sum_{m = k-i}^{k}\tlop(a_{m})\heul^{k-m} = 0 \mod \heul^{i+1}$. Define a multiplication, $\square$, on $\alg[c]$ by setting $\tlop(a \square b) = \tlop(a)\tlop(b)$ for $a,b \in \alg$ and extending $c$-linearly. Together equation \eqref{liftsym} and \eqref{binfr} show that $a\square b = a \cinfstar b$. Precisely, \eqref{liftsym} yields $\hPhi(\tilde{a}_{1})\hPhi(\tilde{a}_{2}) = \hPhi(\widetilde{a_{1}\odot a_{2}}) + \sum_{s =1}^{K}\frac{\binom{K}{s}}{(n+2K -s)_{(s)}}\hPhi(\widetilde{\op_{K-s}(a_{1}, a_{2})})\heul^{s}$, and the star product $a_{1} \square a_{1} = a_{1} \odot a_{2} + \sum_{r \geq 1}B^{\square}_{r}(a_{1}, a_{2})\heul^{r}$ is defined by $\hPhi(\widetilde{a_{1} \square a_{2}}) = \hPhi(\tilde{a}_{1})\hPhi(\tilde{a}_{2})$ (where both $\hPhi$ and the invariant lift have been extended $\heul$-linearly). Equation \eqref{binfr} shows that $B_{r}^{\square}(a_{1}, a_{2}) = B_{r}^{\infty}(a_{1}, a_{2})$, and so $\square = \cinfstar$. This shows that $\cinfstar$ is obtained by transporting the multiplication on $\pol_{0}(T^{\ast}\form)$ to $\alg$ via the invariant lift determined by $\en$. This construction makes evident the associativity of $\cinfstar$.
\end{proof}

\begin{example}\label{infmultex}
Let the notation be as in Section \ref{onedimsection}. In \cite{Cohen-Manin-Zagier} a commutative associative multiplication, $\cmzinf(u_{1}, u_{2})$, of $u_{i} \in \emfb[2k_{i}]$, is defined by a limiting procedure and shown to have the form
\begin{align*}
&\cmzinf(u_{1}, u_{2}) = u_{1}u_{2} + \sum_{r \geq 1} t_{2r}^{\infty}\rc_{2r}(u_{1}, u_{2}),&\\
 &t_{2r}^{\infty} = \frac{(-1/2)_{(r)}}{(r - k_{1} - 1/2)_{(r)}(r - k_{2} - 1/2)_{(r)}(2r - k_{1}-k_{2} -3/2)_{(r)}}.
\end{align*}
Here it will be shown that for the special case $k_{1} = k$, $k_{2} = 0$, the product $\cmzinf(u_{1}, u_{2})$ is the real part of $u_{1}\cinfstar u_{2}$ specialized at $c = 2i$. By \eqref{binfr} and the discussion in Section \ref{onedimsection},
\begin{align*}
&u_{1}\cinfstar u_{2} = \\
&u_{1}u_{2} + \sum_{r \geq 1}\frac{\binom{k}{2r}}{(2k+1-2r)_{(2r)}\binom{2k}{2r}}c^{2r}\rc_{2r}(u_{1}, u_{2}) -  \sum_{r \geq 0}\frac{\binom{k}{2r+1}}{(2k - 2r)_{(2r+1)}\binom{2k}{2r+1}}c^{2r+1}\rc_{2r+1}(u_{1}, u_{2}).
\end{align*}
It is easily checked that
\begin{align*}
&\frac{\binom{k}{2r}}{(2k+1-2r)_{(2r)}\binom{2k}{2r}} = \left(4^{2r}(r -k - 1/2)_{(r)}(2r - k-3/2)_{(r)}\right )^{-1} = (2i)^{-2r}t_{2r}^{\infty},
\end{align*}
and this proves that $\re \{(u_{1} \cinfstar u_{2})_{|c = 2i}\} = \cmzinf(u_{1}, u_{2})$.
\end{example}

\subsection{Isomorphisms of Star Algebras}
The pullback operator associated to $\phi \in\diff(M)$ acts on $(\alg, \pois)$ by graded Poisson automorphisms. By functoriality of the Thomas ambient connection the quantization maps associated to $\en$ and $\ben = \phi \cdot \en$ are related by conjugation, $\bar{\hop}(a) = \phi^{\ast} \circ \hop((\phi^{\ast})^{-1}(a)) \circ (\phi^{\ast})^{-1}$. Using this a simple formal computation shows that $\phi^{\ast}(a \hstar b) = \phi^{\ast}(a) \bhstar \phi^{\ast}(b)$, so that the pullback operator gives a geometric isomorphism between the star products associated to $\en$ and $\ben$.

\begin{lemma}
If $X \in \vect(M)$ is an infinitesimal automorphism of the projective structure, $\en$, on $M$, then $\lie_{X}$ is an inner derivation of the star algebra $(\alg[[\ih]], \pois, \hstar)$ associated to $\en$.
\end{lemma}
\begin{proof}
  By definition $X$ is the infinitesimal generator of a one-parameter family, $\phi_{t}$, of projective automorphisms of $\en$, for which $\phi_{t}^{\ast}(a \hstar b) = \phi_{t}^{\ast}(a) \hstar \phi_{t}^{\ast}(b)$. Differentiating this relation implies that $\lie_{X}(a \hstar b) = \lie_{X}(a) \hstar b + a \hstar \lie_{X}(b)$, showing that $\lie_{X}$ acts as a derivation. Likewise, differentiating $\hop(a) = \phi_{t}^{\ast} \circ \hop((\phi_{t}^{-1})^{\ast}(a)) \circ (\phi_{t}^{-1})^{\ast}$ shows the second equality of
\begin{align*}
\hop(\{X, a\}) = \hop(\lie_{X}a) = [\lie_{X}, \hop(a)] = [\hop(X), \hop(a)] = \hop(X \hstar a - a \hstar X).
\end{align*}
This shows that $\ih \{X, a\} = \ih\lie_{X}a = X \hstar a - a \hstar X$ for all $a \in \alg$, which is the meaning of the statement that $\lie_{X}$ is inner. 
\end{proof}

\begin{proposition}
If two projective structures of dimension at least $2$ determine the same star product, $\hstar$, then the projective structures are related by a projective automorphism. 
\end{proposition}
\begin{proof}
Let projective structures $\ben$ and $\en$ with difference tensor, $\Pi_{ij}\,^{k}$, determine the star products $\bhstar$ and $\hstar$. For $a^{i}, b^{i} \in \emfb^{i}$, computation using \eqref{op011} and \eqref{starvector} shows
\begin{align*}
a \bhstar b - a \hstar b = -\tfrac{\ih^{2}(n+1)}{4(n+2)}\left( \tfrac{2}{1-n}(\nabla_{p}\Pi_{ij}\,^{p} + \tfrac{n+1}{2}\Pi_{ip}\,^{q}\Pi_{jq}\,^{p})a^{(i}b^{j)} + \Pi_{ij}\,^{k}(a^{j}\nabla_{k}b^{i} + b^{j}\nabla_{k}a^{i})\right).
\end{align*}
Suppose that for all vector fields, $a^{i}$, $b^{j}$, there holds $a\bhstar b = a\hstar b$. Fix a point and at that point let the $1$-jet of $b$ be such that $b^{i} = 0$ and $\nabla_{i}b^{j}= c_{i}\,^{j}$. Then at that point $0 = a \Bar{\hstar} b - a \hstar b = -\tfrac{\ih^{2}(n+1)}{4(n+2)}\Pi_{ij}\,^{k}a^{j}c_{k}\,^{i}$. Since this can be made true at every point for all choices of $a^{i}$ and $c_{i}\,^{j}$, it shows that $\Pi_{ij}\,^{k} = 0$, so that the projective structures are the same.
\end{proof}

\begin{theorem}\label{gaugetheorem}
The gauge-equivalence class of the star product $\hstar$ on $(\alg, \{\,,\,\})$ associated to the projective structure, $\en$, does not depend on the choice of projective structure. 
\end{theorem}
\begin{proof}
Given projective structures, $\en$ and $[\bnabla]$, define $D:\alg \to \alg$ by $D = \bar{\hop}^{-1} \circ \hop$ and extend $D$ $\ih$-linearly to $\alg[[\ih]]$. That $D(a \hstar b) = D(a) \bhstar D(b)$ is immediate. For $a \in \emfb^{(i_{1}\dots i_{k})}$, since the principal symbol of each of $\bar{\hop}(a)$ and $\hop(a)$ is $\ih^{k}a$, there holds $D(a) = a + O(\ih)$. Let $D_{1}(a) \in \emfb^{(i_{1}\dots i_{k-1})}$ be the principal symbol of $\ih^{-k}(\bar{\hop}(a) - \hop(a))$, which is a differential operator in $a$ of order at most $k-1$. Thus $D(a) = a + \ih D_{1}(a) + O(\ih)$. Continuing inductively, one finds $D(a) = \sum_{r = 1}^{k}\ih^{r} D_{r}(a)$ where $D_{r}(a) \in \emfb^{(i_{1}\dots i_{k-r})}$. In general $D = \sum_{r \geq 0}\ih^{r}D_{r}$ and $D$ is graded differential in the sense that the restriction to $A_{k}$ of $D_{r}$ is a linear differential operator of order at most $k - r$.\end{proof}

By a Poisson algebra automorphism of $(\alg, \pois)$ is meant a bijective, $\rea$-linear map which is simultaneously an automorphism of the underlying associative algebra structure of $\alg$ and of the Lie algebra structure given by $\pois$. The action of the differential of $\phi \in \diff(M)$ induces an automorphism of $(\alg, \pois)$ to be denoted $\phi^{\ast}$. Each closed one-form, $\alpha$, on $M$, induces an automorphism of $(\alg, \pois)$ as follows. The one-form $\alpha$ determines a $C^{\infty}(M)$-module map $\alg_{1} \to \alg_{0}$ by $A \to \alpha(A)$, and this map may be extended to be identically $0$ on $\alg_{0}$. Because $\alpha :\alg_{0}\oplus\alg_{1} \to \alg_{0}$ is an $\alg_{0}$-module map, it extends uniquely to an associative algebra derivation $\alpha:\alg \to \alg$. Namely, for $A_{1}, \dots, A_{k} \in \alg_{1}$, define $\alpha(A_{1}\odot \dots \odot A_{k}) = \sum_{s = 1}^{k}\alpha(A_{s})A_{1}\odot \dots \odot \hat{A}_{s}\odot \dots \odot A_{k}$ and extend $\rea$-linearly to all of $\alg$. Because $\alpha$ is closed there holds $\alpha(\{A, B\}) = \{\alpha(A), B\} + \{A, \alpha(B)\}$ for $A, B \in \alg_{1}$, and so $\alpha$ extends to a Poisson algebra derivation on all of $\alg$, which is graded nilpotent in the sense that the restriction to $\alg_{k}$ of this derivation is nilpotent of order $k+1$. As the restriction to any $\alg_{k}$ is a finite sum, it makes sense to define $\exp(\alpha) = Id_{\alg} + \sum_{r \geq 1}\tfrac{\alpha^{r}}{r!}$, where $\alpha^{r}$ denotes $r$-fold composition. 

\begin{proposition}[\cite{Grabowski-Poncin}]\label{poissonform}
A Poisson algebra automorphism of $(\alg, \pois)$ has the form $\phi^{\ast} \circ \exp(\alpha)$ for some $\phi \in \diff(M)$ and some closed one-form, $\alpha$, on $M$.
\end{proposition}

\begin{proof}
This follows from Theorem 6 of \cite{Grabowski-Poncin}, though the terminology of Poisson automorphism is used here slightly differently than in that paper. For convenience, the proof is given here in a condensed form. Let $\Phi$ be an automorphism of the Poisson algebra, $(\alg, \pois)$. First it is claimed that $\Phi$ is necessarily filtration preserving in the sense that $\Phi(\alg^{i}) \subset \alg^{i}$. Because $\{\alg_{i}, \alg_{j}\} \subset \alg_{i+j-1}$, elements of $\alg_{0}$ are lowering in the sense that taking the Poisson bracket with an element of $\alg_{0}$ lowers the degree. In fact, $\alg_{0}$ is characterized as the set of $A \in \alg$ such that for any $B \in \alg$ there exists a $k \geq 1$ so that the $k$-fold iterated Poisson bracket $\{A, \{A, \dots \{A, B\} \dots \}\}$ vanishes. This subset is evidently preserved by the action of $\Phi$, and so $\Phi(\alg_{0}) \subset \alg_{0}$. In general, $\alg^{i+1}$ comprises those $A \in \alg$ such that for every $B \in \alg_{0}$ there holds $\{A, B\} \in \alg^{i}$, and from this the claim follows by induction. Precisely, assume $\Phi(\alg^{i}) \subset \alg^{i}$. Then for $A \in \alg^{i+1}$ and $B \in \alg_{0}$, $\{\Phi(A), \Phi(B)\} = \Phi(\{A, B\})$ and as $\{A, B\} \in \alg^{i}$ the inductive hypothesis, implies $\Phi(\{A, B\}) \in A^{i}$. Since the restriction to $\alg_{0}$ of $\Phi$ is an automorpshim of $\alg_{0}$, this shows that $\Phi(A) \in \alg^{i+1}$. By the preceeding $\Phi(\alg_{0}) \subset \alg_{0}$ and $\Phi^{-1}(\alg_{0}) \subset \alg_{0}$ so $\Phi|_{\alg_{0}}$ is an automorphism of the associative algebra $\alg_{0} = C^{\infty}(M)$. Every such automorphism has the form $f \to f\circ \phi^{-1}$ for some $\phi \in \diff(M)$ (see e.g. \cite{Mrcun}). The composition $\phi^{\ast} \circ \Phi$ restricts to the identity on $\alg_{0}$, so it may be assumed that $\Phi|_{\alg_{0}} = Id_{\alg_{0}}$. It will be shown that $\Phi = \exp(\alpha)$ for some closed one-form, $\alpha$. 

If $\Phi|_{\alg_{0}} = Id_{\alg_{0}}$ then $\Phi(fA) = \Phi(f)\Phi(A) = f \Phi(A)$ for all $f \in \alg_{0}$ and $A \in \alg$, so that $\Phi$ is an $\alg_{0}$-module homomorphism. For $A \in \alg_{1}$, $\Phi(A) = \Phi_{1}(A) + \Phi_{0}(A)$ with $\Phi_{i}(A) \in \alg_{i}$, and so $\lie_{A}f = \Phi(\lie_{A}f) = \Phi(\{A, f\}) = \{\Phi(A), \Phi(f)\} = \{\Phi_{1}(A), f\} = \lie_{\Phi_{1}(A)}f$ for all $f \in \alg_{0}$ and $A \in \alg_{1} = \vect(M)$. This implies that $\Phi_{1}(A) = A$ for all $A \in \alg_{1}$. For $A, B \in \alg_{1}$, expanding $\{\Phi(A), \Phi(B)\} = \Phi(\{A, B\})$ yields $\Phi_{0}(\{A, B\})= \{A, \Phi_{0}(B)\} - \{B, \Phi_{0}(A)\}$. It follows that $\Phi_{0}:\alg_{1} \to \alg_{0}$ is an $\alg_{0}$-module homomorphism and a $1$-cocycle in the Lie algebra cohomology of $\vect(M) = \alg_{1}$ with coefficients in $C^{\infty}(M) = \alg_{0}$. Every cocyle of $\vect(M)$ with coefficients in $C^{\infty}(M)$ has the form $c \sdiv + \al$ for some $c \in \rea$ and some closed one-form $\al$ on $M$, where $\sdiv$ denotes any choice of divergence operator (any two divergence operators determine cohomologous cocycles). Evidently such a cocycle is a $C^{\infty}(M)$-module homomorphism if and only if $c = 0$. It follows that $\Phi_{0}(A) = \alpha(A)$ for some closed $1$-form $\alpha$. Since $\alg^{1}$ generates $\alg$ as a $\alg_{0}$-module, two $\alg_{0}$-module homomorphisms of $\alg$ which agree on $\alg^{1}$ must be identical, and so $\Phi = \exp(\alpha)$.
\end{proof}

\begin{lemma}\label{secondorder}
A Poisson algebra automorphism of $(\alg, \pois)$ extends to an automorphism of the star algebra $(\alg[[\ih]], \hstar)$ if and only if it is induced by the action of some element of $\diff(M)$.
\end{lemma}

\begin{proof}
Any $\Phi \in \aut(\alg, \pois)$ acts as an automorphism of $(\alg[[\ih]], \hstar)$ to first order in the sense that $\Phi(a \hstar b) = \Phi(a) \hstar \Phi(b) + O(\ih^{2})$, so to prove the claim it is necessary to examine terms of second order. By Proposition \ref{poissonform} it suffices to show that if $\Phi = \exp(\alpha)$ acts as an automorphism of $(\alg[[\ih]], \hstar)$, then $\alpha = 0$. Let $a^{ij} \in \emfb^{ij}$ and $f \in \emfb[0]$. By \eqref{starvector2} and because $\exp(\alpha)$ acts as the identity on functions on $M$ and as a Poisson algebra automorphism, $\exp(\alpha)(a \hstar f) = \exp(\alpha)(a) \hstar \exp(\alpha)(f)$ if and only if 
\begin{align*}
0 = \op_{0}(\exp(\alpha)(a), f) - \op_{0}(a, f) = \alpha_{j}a^{ij}\nabla_{i}f
\end{align*}
holds for all $a^{ij}$ and all $f$. This implies $\alpha = 0$.
\end{proof}

\begin{lemma}\label{starisogauge}
If $T = \sum_{r \geq 0}\ih^{r}T_{r}$ is a star algebra isomorphism of $(\alg[[\ih]], \hstar)$, then $T_{0}:\alg \to \alg$ is a Poisson algebra automorphism of $(\alg, \pois)$.
\end{lemma}

\begin{proof}
For $a, b \in \alg$, expanding $T(a\hstar b) = T(a)\bhstar T(b)$ by $\ih$-linearity gives 
\begin{align*}
&T(a \odot b) + \ih T(\{a, b\}) =  T(a) \odot T(b) + \ih\{T(a), T(b)\} + O(\ih^{2}).
\end{align*}
Projecting mod $\ih$ this implies $T_{0}(a \odot b) = T_{0}(a) \odot T_{0}(b)$, and then projecting mod $\ih^{2}$ implies that $T_{0}(\{a, b\}) = \{T_{0}(a), T_{0}(b)\}$, so that $T_{0}:\alg \to \alg$ is a Poisson algebra homomorphism. By assumption $T$ has an inverse, $T^{-1} = \sum_{r \geq 0}\ih^{r}S_{r}$. Expanding by $\ih$-linearity gives $a = S_{0}(T_{0}(a)) + O(\ih)$, which implies $S_{0}(T_{0}(a)) = a$, so that $T_{0}$ has an inverse which is a Poisson algebra homorphism, and this proves the claim.
\end{proof}
In particular Lemmas \ref{secondorder} and \ref{starisogauge} imply that a geometric isomorphism of $(\alg[[\ih]], \hstar)$ is induced by extending $\ih$-linearly a Poisson algebra automorphism of $(\alg, \{, \})$ and also that this Poisson algebra automorphism must be induced by a diffeomorphism of $M$. This completes the proof of Theorem \ref{sptheorem}.

\bibliographystyle{amsplain}

\def\cprime{$'$} \def\cprime{$'$} \def\cprime{$'$} \def\cprime{$'$}
  \def\cprime{$'$}
\providecommand{\bysame}{\leavevmode\hbox to3em{\hrulefill}\thinspace}
\providecommand{\MR}{\relax\ifhmode\unskip\space\fi MR }
% \MRhref is called by the amsart/book/proc definition of \MR.
\providecommand{\MRhref}[2]{%
  \href{http://www.ams.org/mathscinet-getitem?mr=#1}{#2}
}
\providecommand{\href}[2]{#2}

\end{document}